# Distributed Random Convex Programming
# via Constraints Consensus

L. Carlone[†], V. Srivastava[‡], F. Bullo[‡], and G.C. Calafiore[†]

### Abstract

This paper discusses distributed approaches for the solution of random convex programs (RCP). RCPs are convex optimization problems with a (usually large) number $N$ of randomly extracted constraints; they arise in several applicative areas, especially in the context of decision under uncertainty, see [2], [3]. We here consider a setup in which instances of the random constraints (the *scenario*) are not held by a single centralized processing unit, but are instead distributed among different nodes of a network. Each node "sees" only a small subset of the constraints, and may communicate with neighbors. The objective is to make all nodes converge to the same solution as the centralized RCP problem. To this end, we develop two distributed algorithms that are variants of the constraints consensus algorithm [4], [5]: the active constraints consensus (ACC) algorithm, and the vertex constraints consensus (VCC) algorithm. We show that the ACC algorithm computes the overall optimal solution in finite time, and with almost surely bounded communication at each iteration of the algorithm. The VCC algorithm is instead tailored for the special case in which the constraint functions are convex also w.r.t. the uncertain parameters, and it computes the solution in a number of iterations bounded by the diameter of the communication graph. We further devise a variant of the VCC algorithm, namely *quantized vertex constraints consensus* (qVCC), to cope with the case in which communication bandwidth among processors is bounded. We discuss several applications of the proposed distributed techniques, including estimation, classification, and random model predictive control, and we present a numerical analysis of the performance of the proposed methods. As a complementary numerical result, we show that the parallel computation of the scenario solution using ACC algorithm significantly outperforms its centralized equivalent.

## I. Introduction

Uncertain optimization problems arise in several engineering applications, ranging from system design, production management, to identification and control, manufacturing and finance, see, e.g., [6]. Uncertainty arises due to the presence of imprecisely known parameters in the problem description. For instance, a system design problem may be affected by the uncertainty in the values of some system components, and control problems can be affected by the inexact knowledge of system model and of the disturbances acting on the system. In the case of uncertain convex optimization problems where the uncertainty in the problem description has a stochastic model (e.g., one assumes random uncertain parameters, with some given probability distribution), the random convex programming (RCP) paradigm recently emerged as an effective methodology to compute "probabilistically robust" solutions, see, e.g., [7], [8], [9].

An instance of an RCP problem typically results in a standard convex programming problem with a large number $N$ of constraints. There are two main reasons for which it is interesting to explore distributed methods for solving RCP instances: first, the number $N$ of constraints may be too large for being stored or solved on a single processing unit; second, there exist application endeavors in which the problem description (objective function and constraints) is naturally distributed among different nodes of an

This work was funded by PRIN grant n. 20087W5P2K from the Italian Ministry of University and Research, by NSF Award CPS 1035917 and by ARO Award W911NF-11-1-0092. The third author thanks Giuseppe Notarstefano for insightful discussions about abstract programming. An early presentation of the active constraints consensus algorithm appeared in [1]. Differences between [1] and the current article include the vertex constraints consensus algorithms, the distributed outliers removal algorithm, proofs and discussions related to the proposed approaches, and the applications to distributed estimation and parallel model predictive control.

[†]L. Carlone and G.C. Calafiore are with the Dipartimento di Automatica e Informatica, Politecnico di Torino, Italy. {luca.carlone,giuseppe.calafiore}@polito.it
[‡]V. Srivastava and F. Bullo are with the Center for Control, Dynamical Systems, and Computation, University of California Santa Barbara, USA. {vaibhav, bullo}@engineering.ucsb.edu

interconnected system. This may happen, for instance, when system constraints depend on measurements acquired by different interacting sensors.

In the last decades, the perspective for solving such large-scale or multi-node problems has switched from *centralized* approaches to *distributed* ones. In the former approach, problem data are either resident on a single node, or transmitted by each node to a central computation unit that solves the (global) optimization problem. In distributed approaches, instead, the computation is fractioned among nodes that must reach a consensus on the overall problem solution through local computation and inter-nodal communication. The advantages of the distributed setup are essentially three-fold: (i) distributing the computation burden and the memory allocation among several processors; (ii) reducing communication, avoiding to gather all available data to a central node; (iii) increasing the robustness of the systems with respect to failures of the central computational unit.

Following this distributed optimization philosophy, we here consider a network of agents or processors that has to solve a random convex program in a distributed fashion. Each node in the network knows a subset of the constraints of the overall RCP, and the nodes communicate with each other with the purpose of determining the solution of the overall problem. Our solution methodology relies on each node iteratively exchanging a small set of relevant constraints, and determining the solution to the RCP in finite time. This methodology is in fact a variation of the constraints consensus algorithm proposed in [4], and further developed in [5].

**Related work.** Distributed and parallel optimization has received significant attention in the literature. In earlier works [10], [11], Lagrangian based decomposition techniques are used to develop decentralized algorithms for large scale optimization problems with separable cost functions. In the seminal work [12], Tsitsiklis investigates the parallel computation of the minimum of a smooth convex function under a setup in which each processor has partial knowledge of the global cost function and they exchange the information of the gradients of their local cost functions to compute the global solution. Recently, Nedíc et. al. [13] generalize the setup of [12] to distributed computation and provide results on the convergence rate and errors bounds for unconstrained problems in synchronous networks. In a similar spirit, Zhu et. al. [14] study primal-dual subgradient algorithm for distributed computation of the optimal solution of a constrained convex optimization problem with inequality and equality constraints. Wei et. al. [15] study a distributed Newton method under a setup in which each node has a partial knowledge of the cost function, and the optimization problem has linear global constraints. Boyd et. al. [16] propose a technique based on dual-decomposition that alternates the updates on different components of the optimization variable. In all these approaches, the proposed algorithms converge to the global solution asymptotically.

An alternative approach to distributed optimization [5], [17], [18] is based on following idea: nodes exchange a small set of constraints at each iteration, and converge in finite time to a consensus set of constraints that determines the global solution of the optimization problem. In particular, Notarstefano et. al. [5] propose constraints consensus algorithm for abstract optimization, while Bürger et. al. [17], [18] present a distributed simplex method for solving linear programs. The algorithms studied in this paper belong to the latter class of algorithms that converge in finite time. Particularly, our first algorithm, the active constraint consensus (ACC), is an adaptation to the RCP context of the constraint consensus algorithm in [5]. Both these algorithms work under similar setups, they have similar approach, and they have very similar properties. The main difference between two algorithms is in the computation of the set of constraints to be transmitted at each iteration. This computation for the algorithm in [5] may require to solve a number of convex programs that grows linearly in the number of constraints and sub-exponentially in the dimension of the problem, while the algorithm considered here always requires the solution of only one convex program. This lower local computation comes at the expense of potentially larger communication at each iteration. In particular, the number of constraints exchanged at each iteration may be higher for the ACC algorithm than the constraints consensus algorithm.

**Paper structure and contributions.** In Section II we recall some preliminary concepts on the constraints of convex programs (support constraints, active constraints, etc.). In Section III we introduce the main distributed random convex programming model, and we describe the setup in which the problem has to be solved. The active constraints consensus algorithm is presented and analyzed in Section IV. In the ACC algorithm, each node at each iteration solves a local optimization problem and transmits to its neighbors the constraints that are tight at the solution (i.e., that are satisfied with equality). We show that the ACC algorithm converges to the global solution in finite time, and that it requires almost surely bounded communication at each iteration. We give some numerical evidence of the fact that the ACC algorithm converges in a number of iterations that is linear in the communication graph diameter. We also provide numerical evidence that parallel implementation of the ACC algorithm significantly reduces the computation time over the centralized computation time. As a side result, we show that the ACC algorithm may distributively compute the solution of any convex program, and that it is particularly effective when the dimension of decision variable is small compared with the number of constraints.

For the special case when the constraints of the RCP are convex in the uncertain parameters, we develop the *vertex constraints consensus* (VCC) algorithm, in Section V. In the VCC algorithm, each node at each iteration constructs the convex hull of the uncertain parameters which define local constraints, and transmits its extreme points to the neighbors. We prove that the VCC algorithm converges to the global solution in a number of iterations equal to the diameter of the communication graph. Moreover, we devise a *quantized vertex constraints consensus* (qVCC) algorithm in which each node has a bounded communication bandwidth and does not necessarily transmit all the extreme points of the convex hull at each iteration. We provide theoretical bounds on the number of the iterations required for qVCC algorithm to converge.

Further, we show in Section VI that each of the proposed algorithms can be easily modified so to enable a *distributed constraints removal strategy* that discards outlying constraints, in the spirit of the RCPV (RCP with violated constraints) framework described in [2]. In Section VII we finally present several numerical examples and applications of the proposed algorithms to distributed estimation, distributed classification, and parallel model predictive control. Conclusions are drawn in Section VIII.

## II. Preliminaries on Convex Programs

Consider a generic $d$-dimensional convex program

$$P[C]: \quad \min_{x \in X} \quad a^\top x \qquad \text{subject to} : \\ f_j(x) \leq 0, \qquad \forall j \in C, \tag{1}$$

where $x \in X$ is the *optimization variable*, $X \subset \mathbb{R}^d$ is a compact and convex domain, $a \in \mathbb{R}^d$ is the *objective direction*, $f_j : \mathbb{R}^d \to \mathbb{R}$, $j \in C$, are convex functions defining problem *constraints*, and $C \subset \mathbb{N}$ is a finite set of indices. We denote the solution of problem $P[C]$ by $x^*(C)$, and the corresponding optimal value by $J^*(C)$; we assume by convention that $x^*(C) = \text{NaN}$ and $J^*(C) = \infty$, whenever the problem is infeasible. We now introduce some definitions, in accordance to [2].

*Definition 1 (**Support constraint set**):* The support constraint set, $\text{Sc}(C) \subseteq C$, of problem $P[C]$ is the set of $c \in C$ such that $J^*(C \backslash \{c\}) < J^*(C)$. □

The cardinality of the set of support constraints is upper bounded by $d+1$, and this upper bound reduces to $d$ if the problem is feasible, see Lemma 2.2 and Lemma 2.3 in [2]. We next provide some definitions.

*Definition 2 (**Invariant and irreducible constraint set**):* A constraint set $S \subseteq C$ is said to be *invariant* for problem $P[C]$, if $J^*(S) = J^*(C)$. A constraint set $S \subseteq C$ is said to be *irreducible*, if $S \equiv \text{Sc}(S)$. □

*Definition 3 (**Nondegenerate problems**):* Problem $P[C]$ is said to be *nondegenerate*, when $\text{Sc}(C)$ is invariant. □

*Definition 4 (**Essential constraint sets**):* An invariant constraint set $S \subseteq C$ of minimal cardinality is said to be an *essential* set for problem $P[C]$. The collection of all essential sets of problem $P[C]$ is denoted as $\text{Es}(C)$. □

*Definition 5 (**Constraints in general position**):* Constraints $f_j(x) \leq 0$, $j \in C$, are said to be in *general position* if the index set $\{i \in C : f_i(x) = 0\}$ has cardinality no larger than $d$, for all $x \in X$. In words, the constraints are in general position if no more than $d$ of the $f_j(x) = 0$ surfaces intersect at any point of the domain $X$. $\square$

*Definition 6 (**Active constraint set**):* The active constraint set $\mathtt{Ac}(C) \subseteq C$ of a feasible problem $P[C]$ is the set of constraints that are tight at the optimal solution $x^*(C)$, that is, $\mathtt{Ac}(C) = \{j \in C : f_j(x^*(C)) = 0\}$. By convention, the active constraint set of an infeasible problem is the empty set. $\square$

Feasible convex programs may have more than one solution, i.e., several values of the optimization variable may attain the same optimal objective value. The convex program $P[C]$ satisfies the *unique minimum condition*, if problem $P[C_i]$ admits a unique solution, for any $C_i \subseteq C$. A convex program that does not satisfy unique minimum condition can be modified into an equivalent problem that satisfies the unique minimum condition, by applying a suitable *tie-breaking rule* (e.g., choosing the lexicographic smallest solution within the set of optimal solutions), see [2]. Accordingly and without loss of generality, in the following we consider convex programs satisfying the unique minimum condition.

### A. Properties of the constraint sets

We now study some properties of the constraint sets in a convex program. We first state the properties of monotonicity and locality in convex programs and then establish some properties of the constraint sets.

*Proposition 1 (**Monotonicity & Locality**, [19], [2]):* For the convex optimization problem $P[C]$, constraint sets $C_1, C_2 \subseteq C$, and a generic constraint $c \in C$, the following properties hold:

i) *Monotonicity:* $J^*(C_1) \leq J^*(C_1 \cup C_2)$;
ii) *Locality:* if $J^*(C_1) = J^*(C_1 \cup C_2)$, then

$$J^*(C_1 \cup \{c\}) > J^*(C_1) \iff J^*(C_1 \cup C_2 \cup \{c\}) > J^*(C_1 \cup C_2). \tag{2}$$

Let the number of different essential sets in $C$ be $n_e$ and $\mathtt{Es}_i(C)$ be the $i$-th essential set. We now state the following proposition on the relationships between support, essential, and active constraint sets.

*Proposition 2 (**Properties of the constraint sets**):* The following statements hold for the constraint sets of a feasible problem $P[C]$:

i) The set of active constraints contains the set of support constraints, that is, $\mathtt{Ac}(C) \supseteq \mathtt{Sc}(C)$;
ii) The set of active constraints contains the union of all the essential sets, that is, $\mathtt{Ac}(C) \supseteq \cup_{i=1}^{n_e} \mathtt{Es}_i(C)$;

**Proof.** See Appendix A.1.
We now state an immediate consequence on Proposition 2.

*Corollary 1 (**Invariance of active constraint set**):* The active constraint set of problem $P[C]$ is an invariant constraint set for $P[C]$.

**Proof.** The second statement of Proposition 2 guarantees that, for any essential set $\mathtt{Es}_i(C)$ of problem $P[C]$, it holds $\mathtt{Ac}(C) \supseteq \mathtt{Es}_i(C)$. By monotonicity, the previous relation implies that (i) $J^*(\mathtt{Ac}(C)) \geq J^*(\mathtt{Es}_i(C))$. However, by the definition of essential set and by monotonicity, we obtain (ii) $J^*(\mathtt{Es}_i(C)) = J^*(C) \geq J^*(\mathtt{Ac}(C))$. Combining (i) and (ii) we prove that $J^*(\mathtt{Ac}(C)) = J^*(C)$, hence the set $\mathtt{Ac}(C)$ is an invariant constraint set for $P[C]$. $\square$

## III. Distributed Random Convex Programming

In this section, we first recall some basic concepts on (standard) random convex programming, [2], and then we define our setup for distributed random convex programming in Section III-B.

### A. Definition and properties of RCPs

A random convex program is a convex optimization problem of the form

$$P[C]: \quad \min_{x \in X} \quad a^\top x \quad \text{subject to :}$$
$$f(x, \delta^{(j)}) \leq 0, \quad j \in C \doteq \{1, \ldots, N\}, \tag{3}$$

where $\delta^{(j)}$ are $N$ independent identically distributed (iid) samples of a random parameter $\delta \in \Delta \subseteq \mathbb{R}^\ell$ having probability distribution $\mathbb{P}$, and $f(x, \delta) : \mathbb{R}^d \times \Delta \to \mathbb{R}$ is convex in $x$, for any $\delta \in \Delta$ (the dependence of $f$ on $\delta$ can instead be generic). The multi-sample $\omega \doteq \{\delta^{(1)}, \delta^{(2)}, \ldots, \delta^{(N)}\}$ is called a *scenario*, and the solution of problem (3) is called a *scenario solution*. Notice that, for given $\omega$, an instance of the RCP (3) has precisely the format of the convex program in (1), for $f_j(x) \doteq f(x, \delta^{(j)})$, and for this reason, with slight abuse of notation, we kept the name $P[C]$, for (3).

A key feature of a RCP is that we can bound a priori the probability that the scenario solution remains optimal for a further realization of the uncertainty [2]. We introduce the following definition.

*Definition 7 (**Violation probability**):* The *violation probability* $V^*(\omega)$ of the RCP (3) is defined as

$$V^*(\omega) \doteq \mathbb{P}\{\delta \in \Delta : \ J^*(\omega \cup \{\delta\}) > J^*(\omega)\},$$

where, $J^*(\omega)$ denotes the optimal value of (3), and $J^*(\omega \cup \{\delta\})$ denotes the optimal value of a modification of problem (3), where a further random constraint $f(x, \delta) \leq 0$ is added to the problem. $\qquad \square$

If problem (3) is nondegenerate with probability one, then the violation probability of the solution satisfies

$$\mathbb{P}\{\omega \in \Delta^N : \ V^*(\omega) \leq \epsilon\} \geq 1 - \Phi(\epsilon; \zeta - 1, N), \tag{4}$$

where $\Phi(\epsilon; q, N) \doteq \sum_{j=0}^{q} \binom{N}{j} \epsilon^j (1 - \epsilon)^{N-j}$ is the cumulative distribution function of a binomial random variable, and $\zeta$ is equal to $d$, if the problem is feasible with probability one, and is equal to $d+1$, otherwise; see Theorem 3.3 of [2]. Furthermore, if one knows a priori that problem (3) is feasible with probability one, then the violation probability $V^*(\omega)$ also represents the probability with which the optimal solution $x^*(\omega)$ of (3) violates a further random constraint, that is

$$V^*(\omega) = \mathbb{P}\{\delta \in \Delta : \ f(x^*(\omega), \delta) > 0\},$$

see Section 3.3 in [2].

For a given $\beta \in (0, 1)$, the bound in equation (4) is implied by

$$\mathbb{P}\{\omega \in \Delta^N : \ V^*(\omega) \leq 2(\log \beta^{-1} + \zeta - 1)/N\} \geq 1 - \beta. \tag{5}$$

In practice, one chooses a *confidence level* $1 - \beta$ close to 1 and picks $N$ large enough to achieve a desired bound on the probability of violation. These bounds on the violation probability neither depend on the uncertainty set $\Delta$, nor on the probability distribution of $\delta$ over $\Delta$. Hence, the RCP framework relaxes basic assumptions underlying robust and chance-constrained optimization [2].

### B. A distributed setup for RCPs

We next describe a distributed formulation of an RCP problem instance. The proposed formulation is similar to the distributed abstract optimization setup in [4], [5]. Consider a system composed of $n$ interacting nodes (e.g., *processors*, *sensors* or, more generically, *agents*). We model inter-nodal communication by a directed graph $\mathcal{G}$ with vertex set $\{1, \ldots, n\}$: a directed edge $(i, j)$ exists in the graph if node $i$ can transmit information to node $j$. We assume that the directed graph $\mathcal{G}$ is strongly connected, that is, it contains a directed path from each vertex to any other vertex. Let $\mathcal{N}_{\text{in}}(i)$ and $\mathcal{N}_{\text{out}}(i)$ be the set of incoming and outgoing neighbors of agent $i$, respectively. Let the *diameter* of the graph $\mathcal{G}$ be $\texttt{diam}(\mathcal{G})$. We state the distributed random programming problem as follows:

*Problem 1 (**Distributed random convex programming**):* A networked system with a strongly connected communication graph has to compute the scenario solution for the random convex program (3), under the following setup:

i) each node knows the objective direction $a$;

ii) each node initially knows only a subset $C_i \subset C$ of the constraints of problem (3) *(the local constraint set)*, $\cup_{i=1}^n C_i = C$;

iii) a generic node $i$ can receive information from the incoming neighbors $\mathcal{N}_{in}(i)$ and can transmit information to the outgoing neighbors $\mathcal{N}_{out}(i)$.

Let $N_i \doteq |C_i|$, for each $i \in \{1, \ldots, n\}$, and let $N = |C|$. Since each node only has partial knowledge of problem constraints, it needs to cooperate with the other nodes to compute the solution of $P[C]$. We say that an iteration at a node has initiated, if the node has received the local information from its neighbors. In the following, we assume that, at any iteration $t \in \mathbb{Z}_{\geq 0}$, node $i$ in the network is able to solve local convex optimization problems of the form:

$$P[L_i(t)] : \min_{x \in X} \quad a^\top x \qquad \text{subject to:}$$
$$f_j(x) \leq 0 \qquad \forall j \in L_i(t) \tag{6}$$

where $L_i : \mathbb{Z}_{\geq 0} \to \text{pow}(C)$ is the subset of constraints that is locally known at node $i$ at time $t$ (possibly with $|L_i| \ll |C|$), and $\text{pow}(C)$ represents the set of all subsets of $C$. We refer to the solution of problem (6) as *local solution* $x_i^*(t) \doteq x^*(L_i(t))$, and the associated value of the objective function as *local optimal value* $J_i^*(t) \doteq J^*(L_i(t))$ (under the convention that $x_i^*(t) = \text{NaN}$ and $J_i^*(t) = \infty$ for infeasible problems).

## IV. Active Constraints Consensus (ACC) algorithm

In this section we describe the active constraints consensus distributed algorithm for solving the RCP (3). We assume that a generic node $i$ at time $t$ can store a small *candidate* constraint set $A_i(t)$, the local optimal solution $x_i^*(t)$, and the local optimal objective $J_i^*(t)$. In the ACC algorithm, each node initially solves the local convex program $P[C_i]$, finds the active constraints $\text{Ac}(C_i)$, and initializes $A_i(0) = \text{Ac}(C_i)$, $x_i^*(0) = x^*(C_i)$, and $J_i^*(0) = J^*(C_i)$. At each iteration $t$ of the algorithm, node $i$ receives the objective values $J_j^*(t)$ and the candidate sets $A_j(t)$ from the incoming neighbors, $j \in \mathcal{N}_{in}(i)$, and builds the constraint set:

$$L_i(t+1) = A_i(t) \cup \left( \cup_{j \in \mathcal{N}_{in}(i)} A_j(t) \right) \cup C_i.$$

Each node then solves problem $P[L_i(t+1)]$ and updates the local quantities, setting $A_i(t+1) = \text{Ac}(L_i(t+1))$, $x_i^*(t+1) = x^*(L_i(t+1))$, and $J_i^*(t+1) = J^*(L_i(t+1))$. The algorithm is iterated until a stopping condition is met (see Remark 1). The details of the algorithm to be executed by each node $i$, $i \in \{1, \ldots, n\}$, are reported as a pseudo code in Algorithm 1. The key properties of the ACC algorithm are summarized in the following proposition.

*Proposition 3 (**Properties of ACC algorithm**):* For a distributed RCP (Problem 1) and the ACC algorithm (Algorithm 1), the following statements hold:

i) the local optimal objective $J_i^*(t)$ is monotonically non-decreasing in the iterations $t$;

ii) the local optimal objective and the local solution converge in a finite number $T$ of iterations to the optimal value $J^*(C)$ and the scenario solution $x^*(C)$ of the RCP;

iii) for each candidate $i$, the local candidate set $A_i(T)$ coincides with the active set $\text{Ac}(C)$ of the RCP;

iv) if constraints in $C$ are in general position, at each iteration of Algorithm 1, each node transmits to each of the outgoing neighbors its current objective value $J_i^*(t)$ and at most $d$ constraints.

**Proof.** The proof of the proposition is an adaptation of the proof of Theorem IV.4 in [5]. We report the proof in Appendix A.2. The main difference in the proofs is that we tailor the demonstration to the exchange of active constraints (instead of the constraints in the *basis*) and we consider explicitly the case of infeasible programs.

*Remark 1 (**Stopping rule for ACC**):* An important fact for the demonstration of claim (i) of Proposition 3 is that if the local optimal objective $J_i^*(t)$ at one node does not change for $2\text{diam}(\mathcal{G}) + 1$ iterations, then convergence has been reached. This fact can be used for implementing a local stopping condition: node $i$ stores an integer ($\text{ncg}$ in Algorithm 1) that counts the number of iterations in which the local optimal objective has not changed. Then the node can stop the algorithm as soon as this counter reaches the value $2\text{diam}(\mathcal{G}) + 1$. The node can also stop iterating the algorithm when an infeasible instance is discovered in its local problem or within the local problems of its neighbors. In particular, as soon a

---

**Algorithm 1:** Active Constraints Consensus (ACC)

---

**Input** : $a$, $C_i$, and $\mathtt{dm} = \mathtt{diam}(\mathcal{G})$;
**Output** : $x^*(C)$, $J^*(C)$, and $\mathtt{Ac}(C)$;

**% Initialization:**
$A_i(0) = \mathtt{Ac}(C_i)$, $\quad J_i^*(0) = J^*(C_i)$, $\quad x_i^*(0) = x^*(C_i)$, $\quad$ and $\mathtt{ncg} = 1$;
$t = 0$;

**% ACC iterations:**
**while** $\mathtt{ncg} < 2\mathtt{dm} + 1$ **and** $J_i^*(t) < \infty$ **do**
    **% Poll neighbors and build:**
    $L_i(t+1) = A_i(t) \cup \left( \cup_{j \in \mathcal{N}_{\mathrm{in}}(i)} A_j(t) \right) \cup C_i$;
    $\bar{J}_i^*(t+1) = \max_{j \in \mathcal{N}_{\mathrm{in}}(i)} J_j^*(t)$;
    **% Check infeasibility:**
    **if** $\bar{J}_i^*(t+1) = \infty$ **then**
        $A_i(t+1) = \emptyset$, $\quad J_i^*(t+1) = \infty$, $\quad x_i^*(t+1) = \mathtt{NaN}$;
        exit;
    **% Update candidate set:**
    $A_i(t+1) = \mathtt{Ac}(L_i(t+1))$, $\quad J_i^*(t+1) = J^*(L_i(t+1))$, $\quad x_i^*(t+1) = x^*(L_i(t+1))$;
    **% Update** $\mathtt{ncg}$ **for stopping condition:**
    **if** $J_i^*(t+1) = J_i^*(t)$ **then**
        $\mathtt{ncg} = \mathtt{ncg} + 1$;
    **else**
        $\mathtt{ncg} = 1$;
    $t = t + 1$;

**return** $x_i^*(t)$, $J_i^*(t)$, $A_i(t)$;

---

node $i$ discovers infeasibility, it sets its objective to $J_i^* = \infty$ and propagates it to the neighbors; as a consequence, all nodes are acknowledged of the infeasibility in at most $\mathrm{diam}(\mathcal{G})$ iterations. $\qquad \square$

*Remark 2 (**Comparison with constraints consensus algorithm [5]**):* The constraint consensus algorithm [5] also distributively computes of the solution of a convex program, and is, in fact, identical to the ACC algorithm whenever the active constraints set and the essential constraints set (basis) are identical. However, in general, the constraint consensus algorithm requires the nodes to compute a basis of the local set of constraints at each iteration, and such computation may be expensive. In particular, for the computation of a basis of a degenerate $d$-dimensional problem with $N_i$ constraints, the algorithm proposed in [5] requires the solution of a number of convex optimization problems that depends linearly on $N_i$ and sub-exponentially on $d$. On the other hand, the active set computation in the ACC algorithm requires the solution of at most one convex program. Particularly, if the local solution $x_i^*(t)$ satisfies all incoming neighbors constraints, then no optimization problem is solved, and the update rule of the ACC algorithm only requires to check if some of the incoming constraints are active. This lower computational expense is achieved at a potentially higher communication. In particular, the ACC algorithm transmits the set of active constraints at each iteration, and the active constraints set is a superset of each basis. $\qquad \square$

*Remark 3 (**Distributed convex programming and constraints exchange**):* The active constraints consensus algorithm can be used for the distributed computation of the solution of any convex program. The distributed strategy is particularly advantageous when the dimension of the decision variable is small and the number of constraints is large (as in the RCP setup), since in this case the nodes only exchange a small subset of constraints of the local constraint sets. Moreover, each constraint $f_j(x) \doteq f(x, \delta^{(j)})$ of an RCP is parameterized in the realization $\delta^{(j)}$, therefore "exchanging" the constraint $f_j(x)$ reduces to transmitting the vector $\delta^{(j)} \in \mathbb{R}^\ell$. $\qquad \square$

## V. Vertex Constraints Consensus (VCC) Algorithms

In this section, we propose distributed algorithms for RCPs, specialized to the case of constraints that are convex in the parameter $\delta$.

*Assumption 1 (**Convex uncertainty**):* For any given $x \in X$, the function $f(x, \delta)$ in (3) is convex in $\delta \in \Delta$. $\qquad\square$

Consider the random convex program in equation (3). Let the *feasible set* of problem $P[C]$ be $\mathtt{Sat}(C) \doteq \{x \in X : f(x, \delta^{(j)}) \leq 0, \forall \ j \in C\}$. Let $\mathtt{co}(C)$ denote the convex hull of uncertainty vectors $\delta^{(j)} \in \Delta$, $j \in C$, and let $\mathtt{vert}(C) \subseteq C$ denote the indices of the uncertainty vectors that form the vertices of $\mathtt{co}(C)$. The following fact, which is a direct consequence of the Jensen's inequality for convex functions, holds.

*Fact 1 (**Invariance of the vertex set**):* If problem $P[C]$ in (3) satisfies Assumption 1, then $\mathtt{vert}(C) \subseteq C$ is an invariant constraint set. $\qquad\square$

As a consequence of the above fact, solving problem $P[\mathtt{vert}(C_i)]$ is equivalent to solving problem $P[C]$. We now present the VCC algorithm.

### A. The VCC algorithm

The VCC algorithm assumes that at time $t$ a generic node $i$ in the network can store a *candidate set* $V_i(t)$, which is initialized to $V_i(0) = \mathtt{vert}(C_i)$ (i.e., it computes the convex hull of the vectors $\delta^{(j)}$, $j \in C_i$, and stores the indices of the vectors being vertices of the convex hull). At each iteration $t$ of the VCC algorithm, node $i$ receives the candidate sets $V_j(t)$ from the incoming neighbors, $j \in \mathcal{N}_{\text{in}}(i)$, and builds the constraint set $L_i(t+1) = V_i(t) \cup \left( \cup_{j \in \mathcal{N}_{\text{in}}(i)} V_j(t) \right)$. Then, the node updates its candidate set with the following rule: $V_i(t+1) = \mathtt{vert}\left( L_i(t+1) \right)$. The algorithm is iterated for $\mathtt{diam}(\mathcal{G})$ iterations, as summarized in Algorithm 2.

---

**Algorithm 2:** Vertex Constraints Consensus (VCC)

---

**Input** : $a$, $C_i$, and $\mathtt{dm} = \mathtt{diam}(\mathcal{G})$;
**Output** : $x^*(C)$, $J^*(C)$, and $\mathtt{vert}(C)$;

% **Initialization:**
$V_i(0) = \mathtt{vert}(C_i)$;
$t = 0$;

% **VCC iterations:**
**while** $t < \mathtt{dm}$ **do**
    % **Poll neighbors and build:**
    $L_i(t+1) = V_i(t) \cup \left( \cup_{j \in \mathcal{N}_{\text{in}}(i)} V_j(t) \right)$;
    % **Update candidate set:**
    $V_i(t+1) = \mathtt{vert}\left( L_i(t+1) \right)$;
    $t = t + 1$;
% **Compute optimal solution and optimal objective:**
$x_i^*(t) = x^*(V_i(t)), \quad J_i^*(t) = J^*(V_i(t))$;

**return** $x_i^*(t)$, $J_i^*(t)$, $V_i(t)$;

---

*Proposition 4 (**Properties of the VCC algorithm**):* For a distributed random convex program (Problem 1) that satisfies Assumption 1, and the VCC algorithm (Algorithm 2), the following statements hold:

i) the local optimal objective $J_i^*(t) \doteq J^*(V_i(t))$ is monotonically non-decreasing in the iterations $t$;

ii) in $T \leq \mathtt{diam}(\mathcal{G})$ iterations the local solution at a generic node $i$ coincides with the scenario solution of the RCP;

iii) for each node $i$ the local candidate set $V_i(T)$ satisfies $V_i(T) = \mathtt{vert}(C) \supseteq \mathtt{Sc}(C)$.

**Proof.** See Appendix A.3.

*Remark 4 (**Computational complexity of convex hull**):* At each iteration of the VCC algorithm each node computes and transmits the convex hull of a set of vectors in $\mathbb{R}^\ell$. There is an extensive literature on the complexity of convex hull computation and on the expected number of vertices in the convex hull, see, e.g., [20], [21], [22]. In particular, it is known that the convex hull of $N$ points in $\mathbb{R}^\ell$ can be computed in $\mathcal{O}(N \log N + N^{\lceil \ell/2 \rceil})$ iterations. Moreover, there exists a $\mathcal{O}(N)$ deterministic algorithm (see [22]) for computing the convex hull of $N$ points uniformly sampled from the interior of a $\ell$-dimensional polytope, and this convex hull has $\mathcal{O}((\log N)^{\ell-1})$ expected number of vertices. □

*Remark 5 (**Distributed uncertain linear programs**):* A remarkable context in which the VCC algorithm can be applied is that of uncertain linear programs. Consider an RCP instance of a standard-form uncertain LP

$$
\begin{aligned}
\min_{x \in X} \quad & a^\top x \qquad \text{subject to :} \\
& u_i^\top\big(z^{(j)}\big)\, x \; - \; v_i\big(z^{(j)}\big) \le 0, \text{ for each } i \in \{1, \ldots, r\}, \text{ and } j \in \{1, \ldots, N\},
\end{aligned}
\tag{7}
$$

where $z^{(j)}$ are iid realizations of some random uncertain parameter $z \in \mathcal{Z}$, where $\mathcal{Z}$ is some arbitrary space, entering the data $u_i(z) \in \mathbb{R}^d$, $v_i(z) \in \mathbb{R}$ in an *arbitrary* way. This RCP does not satisfy Assumption 1 in general, since $u_i(z)$, $v_i(z)$ may be generic nonconvex functions of $z$. However, the problem is readily re-parameterized as

$$
\begin{aligned}
\min_{x \in X} \quad & a^\top x \qquad \text{subject to :} \\
& \delta_i^{(j)}[x^\top\ 1]^\top \le 0, \text{ for each } i \in \{1, \ldots, r\}, \text{ and } j \in \{1, \ldots, N\},
\end{aligned}
\tag{8}
$$

where we defined the parameters $\delta_i = \delta_i(z) \doteq [u_i^\top(x)\quad v_i(z)] \in \mathbb{R}^{1 \times d+1}$. Clearly, each constraint $\delta_i^{(j)}[x^\top\ 1]^\top \le 0$ is now a linear function of $\delta_i$, hence Assumption 1 is satisfied, and the VCC algorithm can be applied to problem (8), operating on the vertices of the convex hull of the $\delta_i^{(j)}$ parameters. Also, problem (8) can be formally cast in the standard RCP format of (1) by setting $f(x, \delta) = \max_{i \in \{1, \ldots, r\}} \delta_i[x^\top\ 1]$, where $\delta$ contains the collection of the $\delta_i$, $i \in \{1, \ldots, r\}$. □

*Remark 6 (**Constraints reexamination**):* The ACC algorithm requires each node $i$ to reexamine its local constraint set $C_i$ at each iteration. This reexamination is attributed to the fact that a constraint that is not active at a given iteration may become active at a later iteration (see [5] for a similar argument for constraints consensus algorithm). The VCC algorithm, instead, requires the knowledge of $C_i$ only for the initialization, and utilizes only the current candidate set and new received constraints to determine the new candidate set. At a generic iteration $t$ of the VCC algorithm at node $i$, any constraint that lies in the interior of the computed convex hull $\mathtt{co}(L_i(t))$ will never belong to any candidate set at future iterations, and therefore, it can be discarded. □

We conclude this section by noticing that the update rule of the VCC algorithm is independent on the objective direction $a$. Therefore, each node does not need to know the objective direction to reach consensus on the set of constraints defining the feasible set of problem $P[C]$.

### B. Quantized VCC algorithm

The size of the constraint set to be transmitted at each iteration of the VCC algorithm may grow exponentially with the dimension of the parameter vector. Such communication at each iteration of the algorithm may not be sustainable for nodes with a limited communication bandwidth. In this section, we address this issue and modify the VCC algorithm to develop the *quantized VCC* (qVCC) algorithm. The qVCC algorithm differs from the VCC algorithm on the following fronts: (i) each node can transmit at most a fixed number $m$ of constraints in a single communication round (*bounded communication bandwidth*); and (ii) a generic node $i$ at time $t$ stores an ordered set, called transmission set, $T_i(t)$, along with the *candidate set*, $V_i(t)$. The algorithm works as follows. Each node initializes $V_i(0) = T_i(0) = \mathtt{vert}(C_i)$, i.e., both sets contain the indices of the constraints corresponding to the vertices of the convex hull $\mathtt{co}(C_i)$. At each iteration $t$ of the qVCC algorithm, each node selects the first $m$ constraints in $T_i(t)$, defining the *current message* $M_i(t)$, and transmits $M_i(t)$ to the outgoing neighbors. When a node receives the messages $M_j(t)$ from the incoming neighbors, $j \in \mathcal{N}_{\text{in}}(i)$, it builds the constraint set $L_i(t+1) = V_i(t) \cup \left( \cup_{j \in \mathcal{N}_{\text{in}}(i)} M_j(t) \right)$. Then, node $i$ updates its candidate set with the following rule: $V_i(t+1) = \text{vert}\left( L_i(t+1) \right)$. Moreover, it updates the transmission set with the rule: $T_i(t+1) = T_i(t) \backslash \{ M_i(t) \cup \left( V_i(t) \backslash V_i(t+1) \right) \} \oplus \{ V_i(t+1) \backslash V_i(t) \}$, where $\oplus$ denotes the concatenation of two ordered sets. Roughly speaking, the updated transmission set, $T_i(t+1)$, is obtained from the previous one, $T_i(t)$, by removing (i) the constraints transmitted at time $t$, i.e., $M_i(t)$, (ii) the constraints that disappeared from the candidate set after the update, i.e., $V_i(t) \backslash V_i(t+1)$, and adding the constraints that became part of the candidate set after the update, $V_i(t+1) \backslash V_i(t)$. Note that the set $T_i(t)$ has to be ordered to implement a first-in-first-out (FIFO) strategy for transmitting constraints to the neighbors. The algorithm is iterated until a stopping condition is met (see Corollary 2). The qVCC algorithm for node $i$ is summarized in Algorithm 3.

---

**Algorithm 3:** Quantized Vertex Constraints Consensus (qVCC)

---

**Input** : $a$, $C_i$, $\text{dm} = \text{diam}(\mathcal{G})$, $m$;
**Output** : $x^*(C)$, $J^*(C)$, and $\text{vert}(C)$;

**% Initialization:**
$V_i(0) = \text{vert}(C_i)$, $\quad T_i(0) = \text{vert}(C_i)$, $\quad$ and stop=0;
$t = 0$;

**% qVCC iterations:**
**while** stop $= 0$ **do**
  **% Build local message $M_i(t)$ by selecting the first $m$ constraints in $T_i(t)$**
  **% Poll neighbors and build:**
  $L_i(t+1) = V_i(t) \cup \left( \cup_{j \in \mathcal{N}_{\text{in}}(i)} M_j(t) \right)$;
  **% Update candidate set and transmission set:**
  $V_i(t+1) = \text{vert}(L_i(t+1))$;
  $T_i(t+1) = T_i(t) \backslash \{ M_i(t) \cup \left( V_i(t) \backslash V_i(t+1) \right) \} \oplus \{ V_i(t+1) \backslash V_i(t) \}$;
  **% Check stopping condition:**
  **if** (all nodes have empty transmission set) **then**
    stop $= 1$;
  $t = t + 1$;
**% Compute optimal solution and optimal objective:**
$x_i^*(t) = x^*(V_i(t))$, $\quad J_i^*(t) = J^*(V_i(t))$;
**return** $x_i^*(t)$, $J_i^*(t)$, $V_i(t)$;

---

Properties of the qVCC algorithm are summarized in Proposition 5. Here, we let $N_{\max}$ be the maximum number of local constraints assigned to a node, i.e., $N_{\max} = \max_{i \in \{1,\dots,n\}} N_i$, and let $d_{\max}$ be the maximum in-degree of a node in the network, i.e., $d_{\max} = \max_{i \in \{1,\dots,n\}} |\mathcal{N}_{\text{in}}(i)|$.

*Proposition 5 (**Properties of qVCC algorithm**):* For a distributed random convex program (Problem 1) that satisfies Assumption 1, and the qVCC algorithm (Algorithm 3), the following statements hold:

  i) The local optimal objective function $J_i^*(t) \doteq J^*(V_i(t))$ is monotonically non-decreasing in the iterations $t$;
  ii) In $T \leq \lceil \frac{N_{\max}}{m} \rceil \frac{(d_{\max}+1)^{\text{diam}(\mathcal{G})}-1}{d_{\max}}$ iterations, the local solution at a generic node $i$ converges to the scenario solution of the RCP;
  iii) For each node $i$ the local candidate set $V_i(T)$ satisfies $V_i(T) = \text{vert}(C) \supseteq \text{Sc}(C)$;

  **Proof.** See Appendix A.4.
We notice that the upper bound on $T$ obtained in Proposition 5 corresponds to the worst case in which all constraints in the local sets need to be transmitted among the nodes. In practice, this bound may be pessimistic, then it is of interest to provide a stopping rule that allows nodes to autonomously detect convergence. We now present an example of stopping rule.

*Corollary 2 (**Stopping rule for qVCC**):* For the qVCC algorithm, if at time $t$ all the transmission sets $T_i(t), i \in \{1, \ldots, n\}$ are empty, then the qVCC algorithm has converged to the scenario solution of the random convex program $P[C]$. Moreover, the situation in which the transmission sets of all nodes are empty can be autonomously detected by each node in $\mathtt{diam}(\mathcal{G})$ iterations.

**Proof.** If at time $t$ the transmission sets are empty, a generic node $i$ satisfies $V_i(t + 1) = V_i(t)$ (no message is received from the incoming neighbors). Moreover, the update rule of the transmission set becomes $T_i(t + 1) = T_i(t) \backslash \{M_i(t) \cup (V_i(t) \backslash V_i(t + 1))\} \oplus \{V_i(t + 1) \backslash V_i(t)\} = \emptyset$. Therefore, the local candidate set and the transmission set remain unchanged for all future iterations, i.e., the qVCC algorithm has converged.

Regarding the second statement, we notice that each node having non-empty transmission set can communicate to all other nodes this situation in $\mathtt{diam}(\mathcal{G})$ iterations. Therefore, if for $\mathtt{diam}(\mathcal{G})$ iterations no node notifies that the local transmission set is non-empty, all transmission sets need be empty, and convergence is reached. $\qquad\square$

## VI. Distributed RCP with Violated Constraints

The RCP framework allows to generalize the probabilistic guarantees of the scenario solution to the case in which $r$ constraints are purposely violated with the aim of improving the objective value $J^*(C)$. Given a problem $P[C]$ and a set $R_r \subset C$, with $|R_r| = r$, RCP theory provides a bound for the probability that a future realization of the random constraints violates $x^*(C \backslash R_r)$, see [2].

In this section we study distributed strategies for removing constraints from a random convex program. RCP theory allows generic constraints removal procedures, with the only requirement that the procedure is permutation invariant (i.e., changing the order of the constraints in $C$ must not change the constraints removed by the procedure). We now present a distributed procedure for removing the $r$ constraints. The procedure works as follows: at each outer iteration the nodes perform one of the distributed algorithms presented before (i.e., ACC, VCC, or qVCC). After attaining convergence, each node selects the constraint $c$ with largest Lagrange multiplier (since nodes share the same set of candidate constraints after convergence, they will choose the same constraint), and each node removes the constraint $c$ from the local constraint set. The distributed procedure is then repeated for $r$ outer iterations (i.e., it terminates after removing the desired number of constraints, $r$). The distributed constraints removal procedure is summarized in Algorithm 4. The acronym CC in Algorithm 4 refers to one of the distributed algorithms presented in the previous sections (i.e., ACC, VCC, or qVCC).

---

**Algorithm 4:** *Distributed Constraints Removal*

**Input** : $a$, $C_i$, $\mathtt{dm} = \mathtt{diam}(\mathcal{G})$, and $r$;
**Output** : $x^*(C \backslash R_r)$, $J^*(C \backslash R_r)$, and $R_r$;

**% Initialization:**
$\eta = 0$, $R_\eta = \emptyset$;

**% Outer iterations:**
**while** $\eta < r$ **do**
  **compute** $[x_\eta^*, J_\eta^*, L_\eta] = \mathtt{CC}(a, C_i, \mathtt{dm}, [m])$;
  **select** $c \in L_\eta$ with largest Lagrange multiplier;
  $C_i = C_i \backslash \{c\}$, and $R_{\eta+1} = R_\eta \cup \{c\}$;
  $\eta = \eta + 1$;

**% Compute optimal solution and optimal objective:**
$[x_r^*, J_r^*, L_r] = \mathtt{CC}(a, C_i, \mathtt{dm}, [m])$;

**return** $x_r^*$, $J_r^*$, $R_r$;

---

We now state some properties of distributed constraints removal procedure:

*Proposition 6 (**Distributed constraints removal**):* The distributed constraints removal procedure in Algorithm 4 is permutation invariant. Moreover, if active constraints consensus algorithm is used for distributed computation of the solution to the RCP in Algorithm 4, then the set of removed constraints corresponds to the one computed with the centralized *constraints removal based on marginal costs* [2].

**Proof.** We start by establishing the first statement. We consider the case in which the ACC algorithm is used for implementing the distributed removal procedure. It follows from Proposition 3 that the local candidate set at each node after convergence coincides with the set of active constraints. Both the set of active constraints and the Lagrange multipliers do not depend on the order of the constrains in $C$, therefore the removal procedure is permutation invariant. The permutation invariance of the distributed constraints removal based on the VCC algorithm can be demonstrated using similar arguments. The second statement is a straightforward consequence of the fact that the active constraints are the only ones that have associated Lagrange multipliers greater than zero (*complementary slackness*); therefore, after performing the ACC algorithm, each node is guaranteed to know all the constraints with nonzero Lagrange multipliers, from which it can select the one with largest multiplier. □

We conclude this section with some comments on the trade-off between the use of the ACC and the VCC algorithm in the distributed removal procedure (Algorithm 4). First of all we notice that the ACC algorithm is able to return a constraint set only in feasible problems (otherwise the active constraint set is empty, by convention); therefore, the ACC-based removal procedure applies only to feasible problem instances. On the other hand, under Assumption 1, the VCC-based removal procedure applies in the infeasible case as well. However, when using the VCC (or the qVCC), it is not possible to establish the parallel with the centralized case, since it is possible to have constraints with non-zero Lagrange multipliers that are not in the set computed by the VCC algorithm.

## VII. Applications and Numerical Examples

### A. Distributed ellipsoidal estimation

In this section we discuss the problem of determining a confidence ellipsoid for an unknown random parameter. We study this problem considering three settings: (i) nodes in a network can directly measure the parameter (Section VII-A1), (ii) nodes can measure a linear function of the parameter (Section VII-A2), (iii) nodes may take linear measurements of the parameter using possibly different measurement models (Section VII-A3).

*1) Computing a confidence ellipsoid:* In this section we discuss the problem of determining a confidence ellipsoid for an unknown random parameter $y \in \mathbb{R}^q$ for which $N$ iid realizations $y^{(j)}$, $j \in \{1, \ldots, N\}$ are available. We consider first the case in which all the $N$ realizations are collected at a single unit that solves the problem in a centralized way, and then outline a distributed setup of this problem in Remark 7.

A generic (bounded) ellipsoid, parameterized in its center $\hat{y} \in \mathbb{R}^q$ and shape matrix $W_y \in \mathbb{R}^{q \times q}$, $W_y \succ 0$, is represented as

$$\mathcal{E}_y = \{y \in \mathbb{R}^q : (y - \hat{y})^\top W_y (y - \hat{y}) \leq 1\}. \tag{9}$$

As a measure of size of $\mathcal{E}_y$ we consider the volume, which is proportional to the square root of the determinant of $W_y^{-1}$. Then, the problem of finding the smallest ellipsoid enclosing the given realizations can be stated in the form of the following convex optimization problem

$$\min_{\hat{y}, W_y \succ 0} \texttt{logdet}(W_y^{-1}) \qquad \text{subject to :}$$
$$(y^{(j)} - \hat{y})^\top W_y (y^{(j)} - \hat{y}) \leq 1, \quad \text{for each } j \in \{1, \ldots, N\}. \tag{10}$$

The number of variables in this problem is $q(q+3)/2$, corresponding to $q$ variables describing the center $\hat{y}$, plus $q(q+1)/2$ variables describing the free entries in the symmetric matrix $W_y$. We can convert the optimization problem (10) into an equivalent one having linear cost function by introducing a slack variable (see Remark 3.1 in [2]); the dimension of the problem with linear objective is then $d = q(q+3)/2 + 1$.

Since the realizations $y^{(j)}$ are assumed random and iid, problem (10) clearly belongs to the class of RCPs. Moreover, this problem is always feasible, and its solution is unique (see, for instance, Section 3.3 in [23]). Therefore, we can apply (5) to conclude that with high probability $1 - \beta$ (here, $\beta$ is typically set to a very low value, say $\beta = 10^{-9}$) the ellipsoid computed via (10) is an $(1 - \epsilon)$-confidence ellipsoid for $y$, with $\epsilon = 2(\log \beta^{-1} + d - 1)/N$. In words, we know with practical certainty that $\mathcal{E}_y$ contains $y$ with probability larger than $1 - \epsilon$, i.e., it encloses a probability mass at least $1 - \epsilon$ of $y$. Furthermore, we observe that the constraints in (10) are convex functions also with respect to the "uncertainty" terms $y^{(j)}$, hence this problem satisfies Assumption 1, enabling the application of the VCC or qVCC algorithms.

*Remark 7 (**Distributed computation of measurement ellipsoid**):* The solution to the optimization problem (10) can be computed in distributed fashion using any of the algorithms proposed in this paper, by considering a setup in which $n$ nodes are available, and each node only knows initially $N_i$ local realizations of $y$, with $\sum_{i=1}^{n} N_i = N$. Application of ACC, VCC, or qVCC algorithms entails that each node iteratively exchanges a subset of realizations $y^{(j)}$ with its neighbors in order to reach consensus on the set of realizations defining the optimal solution to (10). □

*2) Ellipsoidal parameter estimation in a linear model:* We now extend the previous setup by considering the case in which linear measurements $y$ of an unknown parameter $\theta$ are used to infer an ellipsoid of confidence for the parameter itself. Consider the classical situation in which $y$ is related to $\theta$ via a linear model

$$y = F\theta, \tag{11}$$

with $F \in \mathbb{R}^{q \times p}$, where $\theta$ is the input parameter, and $y$ is a measured output. Suppose that $\theta^{(1)}, \ldots, \theta^{(N)}$, are $N$ iid realization of the unobservable parameter $\theta$, and that $y^{(1)}, \ldots, y^{(N)}$ are the corresponding observed measurements: $y^{(i)} = F\theta^{(i)}$. We first consider the centralized case, in which a single node uses the measurements to infer an ellipsoid of confidence for $\theta$. Given the observations $y^{(1)}, \ldots, y^{(N)}$, we can compute a unique minimum-size ellipsoid $\mathcal{E}_y$ containing the observations, by solving problem (10). From the reasoning in Section VII-A1 we know with practical certainty that $\mathcal{E}_y$ is a $(1 - \epsilon)$-confidence ellipsoid for $y$. Now, the condition $y \in \mathcal{E}_y$, together with the linear relation in (11), imply that the set of parameters $\theta$ that are compatible with output $y \in \mathcal{E}_y$ is a (possibly unbounded) ellipsoid $\mathcal{E}$ described by the quadratic inequality condition $(F\theta - \hat{y})^\top W_y (F\theta - \hat{y}) \le 1$, that is

$$\begin{bmatrix} \theta \\ 1 \end{bmatrix}^\top \begin{bmatrix} F^\top W_y F & F^\top W_y \hat{y} \\ * & \hat{y}^\top W_y \hat{y} - 1 \end{bmatrix} \begin{bmatrix} \theta \\ 1 \end{bmatrix} \le 0. \tag{12}$$

Since $y \in \mathcal{E}_y$ if and only if $\theta \in \mathcal{E}$, and since with practical certainty $\mathbb{P}\{y \in \mathcal{E}_y\} \ge 1 - \epsilon$, we also have that $\mathbb{P}\{\theta \in \mathcal{E}\} \ge 1 - \epsilon$, hence we found a region $\mathcal{E}$ within which $\theta$ must be contained with probability no smaller than $1 - \epsilon$.

In the next section, we provide an extension of this linear estimation framework to a distributed setup in which $n$ nodes collect linear measurements of $\theta$, using possibly heterogeneous models.

*3) Ellipsoidal parameter estimation in heterogeneous network:* Suppose that there are $n_s$ subsets of nodes, say $\mathcal{V}_1, \ldots, \mathcal{V}_{n_s}$, such that each node in $\mathcal{V}_j$ uses the same linear measurement model

$$y_i = F_j \theta, \quad \text{for each } i \in \mathcal{V}_j, \tag{13}$$

and it collects $N_i$ measurements

$$y_i^{(k)} = F_j \theta^{(k)}, \quad \text{for each } k \in \{1, \ldots, N_i\},$$

where $\theta^{(k)}, k \in \{1, \ldots, N_i\}$, are iid. Moreover, it is assumed that realizations of $\theta$ available at a node $i$ are independent from realizations available at node $j$, for each $i, j$. We here detail the procedure for computing a confidence ellipsoid for $\theta$, by first assuming a centralized case in which all measurements from nodes in $\mathcal{V}_j$ are available at a central node, and then we refer to Remark 8 for outlining the corresponding distributed implementation.

If all measurements from nodes in $\mathcal{V}_j$ are available to a central computational unit, then this unit can first construct (by solving problem (10)) an ellipsoid of confidence $\mathcal{E}_y^j$ for the collective measurements $y_i^{(k)}$, $i \in \mathcal{V}_j$, $k \in \{1, \ldots, N_i\}$:

$$\mathcal{E}_y^j = \{y : (y - \hat{y}_j)^\top W_j (y - \hat{y}_j) \leq 1\},$$

and then infer an ellipsoid of confidence $\mathcal{E}_j$ for $\theta$ according to eq. (12):

$$\mathcal{E}_j = \left\{ \theta \in \mathbb{R}^p : \begin{bmatrix} \theta \\ 1 \end{bmatrix}^\top \begin{bmatrix} F_j^\top W_j F_j & F_j^\top W_j \hat{y}_j \\ * & \hat{y}_j^\top W_j \hat{y}_j - 1 \end{bmatrix} \begin{bmatrix} \theta \\ 1 \end{bmatrix} \leq 0 \right\}.$$

This procedure can be repeated for each $\mathcal{V}_j$, $j \in \{1, \ldots, n_s\}$, thus obtaining $n_s$ ellipsoidal sets $\mathcal{E}_j$ that (with practical certainty) contain $\theta$ with probability no smaller than $1 - \epsilon_j$. "Fusing" the information from all the confidence ellipsoids $\mathcal{E}_j$, a standard probabilistic argument leads to stating that (again with practical certainty) the unknown parameter is contained in the intersection $\mathcal{I} = \cap_{j=1}^{n_s} \mathcal{E}_j$ with probability no smaller than $\mu \doteq \prod_{j=1}^{n_s} (1 - \epsilon_j)$. Clearly, any set that contains the intersection $\mathcal{I}$ has a probability no smaller than $\mu$ of containing $\theta$. We may then find an ellipsoid $\mathcal{E}$ covering the intersection $\mathcal{I}$, as follows. We describe the to-be-computed ellipsoid $\mathcal{E}$ as

$$\begin{bmatrix} \theta \\ 1 \end{bmatrix}^\top \begin{bmatrix} W & W\hat{\theta} \\ * & \hat{\theta}^\top W \hat{\theta} - 1 \end{bmatrix} \begin{bmatrix} \theta \\ 1 \end{bmatrix} \leq 0,$$

where $\hat{\theta}$ is the center of the ellipsoid and $W \succ 0$ is its shape matrix. Then a sufficient condition for $\mathcal{E}$ to contain $\mathcal{I}$ can be obtained through the so-called $\mathcal{S}$-procedure [24]: if there exist $n_s$ scalars $\tau_j \geq 0$, $j \in \{1, \ldots, n_s\}$, such that

$$\begin{bmatrix} W & W\hat{\theta} \\ * & \hat{\theta}^\top W \hat{\theta} - 1 \end{bmatrix} - \sum_{j=1}^{n_s} \tau_j \begin{bmatrix} F_j^\top W_j F_j & F_j^\top W_j \hat{y}_j \\ * & \hat{y}_j^\top W_j \hat{y}_j - 1 \end{bmatrix} \preceq 0,$$

then $\mathcal{E} \supseteq \cap_{j=1}^{n_s} \mathcal{E}_j$. Defining a vector $\tilde{\theta} = W\hat{\theta}$, we can write the previous condition as:

$$\begin{bmatrix} W - \sum_{j=1}^{n_s} \tau_j (F_j^\top W_j F_j) & \tilde{\theta} - \sum_{j=1}^{n_s} \tau_j (F_j^\top W_j \hat{y}_j) \\ * & -1 - \sum_{j=1}^{n} \tau_j (\hat{y}_j^\top W_j \hat{y}_j - 1) \end{bmatrix} + \begin{bmatrix} \mathbf{0}_p \\ \tilde{\theta}^\top \end{bmatrix} W^{-1} \begin{bmatrix} \mathbf{0}_p \\ \tilde{\theta}^\top \end{bmatrix}^\top \preceq 0,$$

where $\mathbf{0}_p$ is a matrix in $\mathbb{R}^{p \times p}$ with all zero entries. Using the Schur complement rule, this latter condition is equivalent to the following LMI in $W$, $\tilde{\theta}$, and $\tau_1, \ldots, \tau_{n_s}$:

$$\begin{bmatrix} W - \sum_{j=1}^{n_s} \tau_j (F_j^\top W_j F_j) & \tilde{\theta} - \sum_{j=1}^{n_s} \tau_j (F_j^\top W_j \hat{y}_j) & \mathbf{0}_p \\ * & -1 - \sum_{j=1}^{n} \tau_j (\hat{y}_j^\top W_j \hat{y}_j - 1) & \tilde{\theta}^\top \\ * & * & W \end{bmatrix} \preceq 0. \tag{14}$$

Then, the shape matrix $W$ of the minimum volume ellipsoid $\mathcal{E} \supseteq \mathcal{I}$ can be computed by solving the following convex program

$$\min_{\hat{\theta}, W \succ 0, \tau_1 \geq 0, \ldots, \tau_{n_s} \geq 0} \texttt{logdet}(W^{-1})$$
$$\text{subject to} : (14). \tag{15}$$

After obtaining the optimal solution of problem (15), the center of the minimum volume ellipsoid can be computed as $\hat{\theta} = W^{-1}\tilde{\theta}$.

*Remark 8 (**Distributed estimation in heterogeneous network**):* A distributed implementation of the procedure previously described goes as follows. We assume that each node $i \in \{1, \ldots, n\}$, knows all the measurement models $\{F_1, \ldots, F_{n_s}\}$, and acquires $N_i$ measurements according to its own model $F_j$, see (13). Each node $i$ then maintains $n_s$ different local constraint sets $C_i^j$, $j \in \{1, \ldots, n_s\}$, simultaneously, and initializes the $j$-th set $C_i^j$ to the local measurements set of node $i$, if $i \in \mathcal{V}_j$, or to the empty set, otherwise. Then, each node runs a distributed constraint consensus algorithm (either ACC, or VCC, or qVCC) simultaneously on each of its local constraint sets. In this way, upon convergence, each node has all the optimal ellipsoids $\mathcal{E}_j$, $j \in \{1, \ldots, n_s\}$. Once this consensus is reached, each node can compute locally the enclosing ellipsoid $\mathcal{E} \supseteq \cap_{j=1}^{n_s} \mathcal{E}_j$, by solving the convex program (15). $\qquad\square$

*4) Numerical results on distributed ellipsoid computation:* We now elucidate on the distributed ellipsoid computation with some numerical examples. In particular, we demonstrate the effectiveness of our algorithms for (i) distributed computation of the enclosing ellipsoid when each node can measure the random parameter $\theta$ with the same measurement model; (ii) parallel computation of the enclosing ellipsoid; and (iii) distributed computation of the enclosing ellipsoid when each node can only measure some components of the random parameter $\theta$.

*Example 1 (**Distributed estimation in homogeneous sensor network**):* Consider the setup in which $n$ sensors measure a random variable $\theta$, using the same measurement model $y = F\theta$ (homogeneous sensor network), where we set for simplicity $F = I_p$ (the identity matrix of size $p$). We assumed $\theta \in \mathbb{R}^2$ to be distributed according to the following mixture distribution:

$$\theta = \begin{cases} \gamma_1 & \text{with probability } 0.95 \\ \gamma_2 + 10\gamma_1 & \text{with probability } 0.05, \end{cases}$$

where $\gamma_1 \in \mathbb{R}^2$ is a standard Normal random vector, and $\gamma_2 \in \mathbb{R}^2$ is uniformly distributed in $[-1, 1]^2$. The overall number of measurements (acquired by all nodes) is $N = 20000$; the size of the local constraint sets is $N/n$. We consider the case in which the nodes in the network solve the RCP in equation (10) using one of the algorithms proposed in this paper. We consider two particular graph topologies: a *chain graph* and a *geometric random graph*. For the geometric random graph, we picked nodes uniformly in the square $[0, 1]^2$ and choose a communication radius $r_c > 2\sqrt{2}\sqrt{\log(n)/n}$ to ensure that the graph is strongly connected with high probability [25]. In Table I we report the maximum number of iterations and the maximum number of exchanged constraints for each algorithm. Statistics are computed over 20 experiments. The ACC algorithm requires nodes to exchange a small number of constraints, and it converges in a number of iterations that grows linearly in the graph diameter. For the VCC algorithm, the maximum number of iterations for convergence is equal to the graph diameter. For the considered problem instances, the number of constraints to be exchanges among the nodes is small. We picked $m = 5$ for the qVCC algorithm. Table I reports the number of iterations required by the qVCC to meet the halting conditions described in Corollary 2.

| | No. Nodes | Diameter | ACC | | VCC | | qVCC | |
|---|---|---|---|---|---|---|---|---|
| | | | Iter. | Constr. | Iter. | Constr. | Iter. | Constr. |
| Geometric random graph | 10 | 1 | 5 | | 1 | | 6 | |
| | 50 | 2 | 7 | | 2 | | 8 | |
| | 100 | 3 | 10 | 6 | 3 | 19 | 9 | 5 |
| | 500 | 5 | 16 | | 5 | | 13 | |
| Chain graph | 10 | 10 | 36 | | 10 | | 21 | |
| | 50 | 50 | 187 | | 50 | | 101 | |
| | 100 | 100 | 375 | 5 | 100 | 23 | 200 | 5 |
| | 500 | 500 | 1910 | | 500 | | 1000 | |

TABLE I
DISTRIBUTED COMPUTATION IN HOMOGENEOUS SENSOR NETWORK: MAXIMUM NUMBER OF ITERATIONS, MAXIMUM NUMBER OF EXCHANGED CONSTRAINTS, AND DIAMETER FOR DIFFERENT GRAPH TOPOLOGIES, AND FOR EACH OF THE PROPOSED ALGORITHMS.

*Example 2 (**Parallel computation of confidence ellipsoid**):* In this example we consider the same setup as in Example 1, but we solve the RCP (10) in distributed fashion assuming a complete communication graph. A complete communication graph describes a parallel computation setup in which each processor can interact with all the others. In this case, we focus on the ACC algorithm. In Fig. 1 we report the dependence of the number of iterations on the number of nodes, number of constraints, and dimension of the parameter $y = \theta$ to be estimated. In the considered problem instances the iterations of the ACC algorithm do not show any dependence on these three factors. In Fig. 2 we show some statistics on the number of exchanged constraints. In particular, we compare the number of constraints exchanged among nodes at each communication round with the dimension $d = p(p+3)/2 + 1$ (recall that $p = q$ in this example) of the RCP (Section VII-A1): in Proposition 3 we concluded that the number of constraints exchanged at each communication round is bounded by $d$. Fig. 2 shows that in the considered problem instances, the number of constraints is below this upper bound, which is shown as a dashed line. For space reasons we do not report results on the dependence of the number of exchanged constraints on the total number of constraints $N$ and on the number of nodes $n$. In our test the number of exchanged constraints was practically independent on these two factors and remained below 5 in all tests.

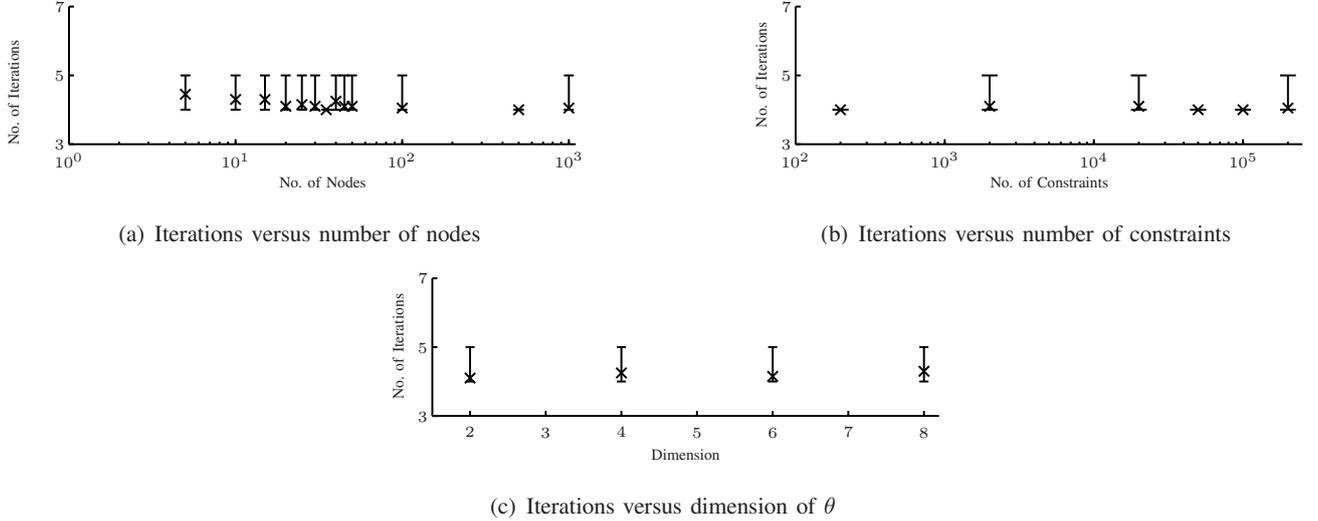

(a) Iterations versus number of nodes

(b) Iterations versus number of constraints

(c) Iterations versus dimension of $\theta$

Fig. 1. Parallel computation of confidence ellipsoid using ACC algorithm. (a) number of iterations required for convergence with different number of nodes, with fixed number of constraints $N = 20000$ and fixed dimension $p = 2$ of $\theta$; (b) number of iterations for different number of constraints, with fixed number of nodes $n = 50$ and fixed dimension $p = 2$; (c) number of iterations for different dimensions of $\theta$, with fixed number of nodes $n = 50$ and number of constraints $N = 20000$. In each figure the cross denotes the average number of iterations, whereas the bar defines the maximum and the minimum observed numbers of iterations.

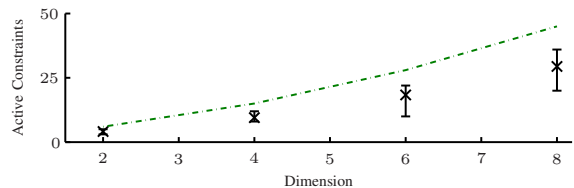

Fig. 2. Parallel computation of confidence ellipsoid using ACC algorithm: (bars) number of constraints exchanged among nodes for different dimension $p$ of $\theta$, with fixed number of constraints $N = 20000$ and fixed number of nodes $n = 50$. The cross denotes the average number of iterations, whereas the bar defines the maximum and the minimum observed numbers of iterations; (dashed line) maximum number of constraints in generic position $d = p(p+3)/2 + 1$.

In Fig. 3 we compare the computational effort required by the ACC algorithm in the parallel setup with a standard centralized solver in charge of solving the convex program (10). We used `CVX/SeDuMi` [26] as a centralized parser/solver, and we compared the computation times required for solving the problem, for different number of nodes, number of constraints, and dimension of the parameter $\theta$. The use of the ACC algorithm provides a remarkable advantage in terms of computational effort. For a large number of constraints, this advantage is significant even for a small number of processors.

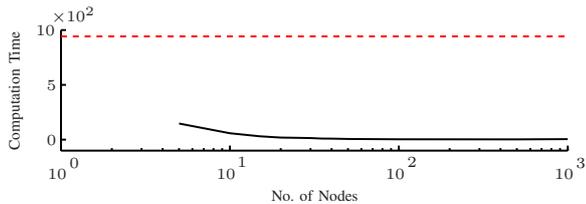

(a) Comp. time versus no. of nodes

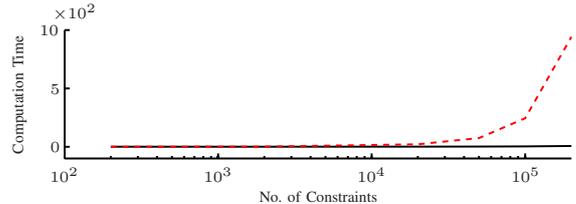

(b) Comp. time versus no. of constraints

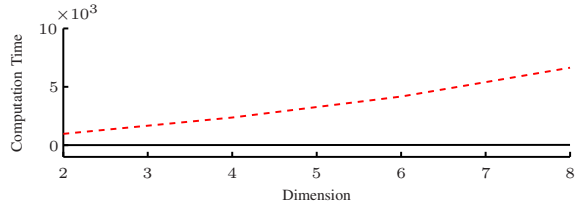

(c) Comp. time versus dimension of $\theta$

Fig. 3. Parallel computation of confidence ellipsoid. The solid black line represents the parallel computation time required for solving the RCP using the ACC algorithm, and dashed red line represents the computation time required for centralized solution of the RCP.

*Example 3 (**Distributed estimation in heterogeneous sensor network**):* We now consider the distributed computation of a parameter ellipsoid in a network with $n$ nodes. We assume that half of the nodes in the network takes measurements of $\theta \in \mathbb{R}^2$ according to the measurement model $y_1 = F_1\theta$, where $F_1 = [1\ 0]$; the remaining nodes use the measurement model $y_2 = F_2\theta$, where $F_2 = [0\ 1]$. We consider $\theta$ distributed according to a mixture distribution, as in Example 1. The nodes acquires 20000 measurements for each measurement model. They then estimate the set $\mathcal{E}$ according to Remark 8. In Table II we report some statistics related to the computation of the sets $\mathcal{E}_1$ and $\mathcal{E}_2$ using the ACC and the VCC algorithms, see Remark 8. After the computation of $\mathcal{E}_1$ and $\mathcal{E}_2$, each node can locally retrieve the set $\mathcal{E}$ solving problem (15), see Fig. 4.

| | No. Nodes | Diameter | ACC | | VCC | |
|---|---|---|---|---|---|---|
| | | | Iter. | Constr. | Iter. | Constr. |
| Geometric random graph | 10 | 1 | 4 | | 1 | |
| | 50 | 2 | 7 | 4 | 2 | 4 |
| | 100 | 3 | 10 | | 3 | |
| | 500 | 5 | 16 | | 5 | |
| Chain graph | 10 | 10 | 28 | | 10 | |
| | 50 | 50 | 148 | 4 | 50 | 4 |
| | 100 | 100 | 298 | | 100 | |
| | 500 | 500 | 1498 | | 500 | |

TABLE II

DISTRIBUTED ESTIMATION IN HETEROGENEOUS SENSOR NETWORK: MAXIMUM NUMBER OF ITERATIONS, MAXIMUM NUMBER OF EXCHANGED CONSTRAINTS, AND DIAMETER FOR DIFFERENT GRAPH TOPOLOGIES, FOR ACC AND VCC ALGORITHMS.

According to Section VII-A3 we can conclude that for $j \in \{1, 2\}$, with confidence level $1-\beta = 1-10^{-8}$, $\mathcal{E}_j$ is a $(1-\epsilon_j)$-confidence ellipsoid for $\theta$, with $\epsilon_j = 2 \cdot 10^{-3}$. Then, with practical certainty the ellipsoid $\mathcal{E}$ is a $\mu$-confidence ellipsoid for $\theta$, with $\mu = (1-\epsilon_1)(1-\epsilon_2) \approx 0.995$.

### B. Distributed linear classification

A classical problem in binary linear classification is to determine a linear decision surface (a hyperplane) separating two clouds of binary labelled multi-dimensional points, so that all points with label $+1$ fall on one side of the hyperplane and all points with label $-1$ on the other side, see Fig. 5. Formally, one is given a set data points (*features*) $b_j \in \mathbb{R}^p$, $j \in \{1, \dots, N\}$, and the corresponding class label $l_j \in \{-1, +1\}$, and seeks a suitable hyperplane $\mathcal{H} = \{s \in \mathbb{R}^p : \theta^\top s + \rho = 0\}$, with $\theta \in \mathbb{R}^p$ and $\rho \in \mathbb{R}$, such that features with different labels belong to different half-spaces w.r.t. $\mathcal{H}$, and the margin of separation between the classes is maximized (maximum margin classifier, see [27]). If the features are linearly separable, then the optimal separating hyperplane solves the following minimization problem [28]:

$$\min_{\theta, \rho} \quad \|\theta\|_2 \qquad \text{subject to:}$$
$$l_j(b_j^\top \theta + \rho) \geq 1, \ \text{for each } j \in \{1, \dots, N\}. \tag{16}$$

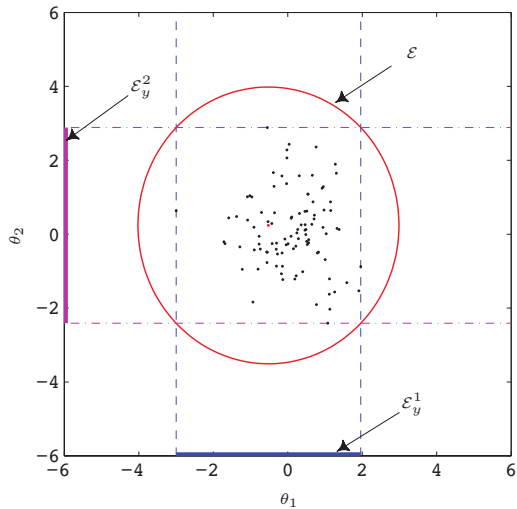

Fig. 4. Distributed estimation in heterogeneous sensor network: the black dots are 100 realizations of the random parameter $\theta = [\theta_1\ \theta_2]^\top$. Nodes with measurement model 1 can measure $y_1 = F_1\theta = [1\ 0]\ \theta = \theta_1$ and compute the corresponding measurement set $\mathcal{E}_y^1$ (shown as a solid blue line), and the set $\mathcal{E}_1$ (the strip delimited by dashed blue lines) of parameters compatible with $\mathcal{E}_y^1$. Similarly, nodes with measurement model 2 can measure $y_2 = F_2\theta = [0\ 1]\ \theta = \theta_2$ from which the network builds the set $\mathcal{E}_y^2$ (shown as a solid magenta line) and the set $\mathcal{E}_2$ (the strip delimited by dashed magenta lines) of parameters compatible with $\mathcal{E}_y^2$. From the sets $\mathcal{E}_1$ and $\mathcal{E}_2$, each node can compute the bounding ellipsoid $\mathcal{E} \supseteq \mathcal{E}_1 \cap \mathcal{E}_2$, by solving problem (15).

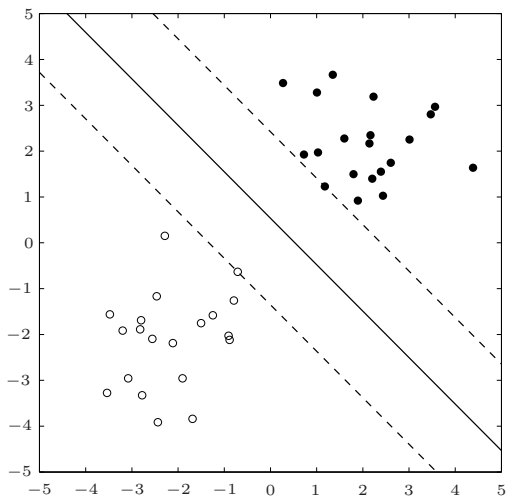

Fig. 5. Binary linear classification: Two clouds of points having labels $+1$ (full circles) and $-1$ (empty circles), respectively, need be separated by a hyperplane $\mathcal{H}$, which maximizes the margin of separation between the classes.

To deal with possibly infeasible problem instances (i.e., non-linearly separable data), it is common to include a slack variable, allowing (but penalizing) misclassification:

$$\min_{\theta, \rho, \nu \geq 0} \quad \|\theta\|_2 + \nu \qquad \text{subject to :}$$
$$l_j({b_j}^\top \theta + \rho) \geq 1 - \nu, \text{ for each } j \in \{1, \ldots, N\}. \tag{17}$$

If the observed datum/label pairs $\delta^{(j)} = (b_j, l_j)$, $j \in \{1, \ldots, N\}$, are interpreted as realization of a random datum/label variable $\delta = (b, l)$, then problem (17) is an instance of the following RCP in dimension $d = p + 3$:

$$\min_{\theta, \rho, \phi \geq 0, \nu \geq 0} \quad \phi \quad \text{subject to :} \tag{18}$$
$$l_j({b_j}^\top \theta + \rho) \geq 1 - \nu, \text{ for each } j \in \{1, \ldots, N\}, \tag{19}$$
$$\|\theta\|_2 + \nu \leq \phi.$$

Such and RCP is always feasible, and it admits a unique optimal solution with probability one, see, e.g., [28]. Therefore, we can apply (5) to conclude that with practical certainty the hyperplane $\mathcal{H}$, obtained as solution of (18), remains an optimal separating hyperplane also after adding a new realization to the training data.

Problem (18) is readily amenable to distributed solution via the ACC algorithm, by assuming that the $N$ constraints in (19) are subdivided into $n$ disjoint subsets of cardinality $N_i$ each, $i \in \{1, \ldots, n\}$, $\sum_{i=1}^n N_i = N$, and that each subset is assigned to a node as local constraint set. The constraints in (19) are linear, hence the problem can also be solved via the VCC or qVCC algorithm, see Remark 5.

*1) Numerical results on distributed linear classification:* We next present numerical examples of distributed and parallel computation of linear classifier.

*Example 4 (**Distributed linear classification**):* In this section we consider the case in which the training set $\delta^{(j)} = (b_j, l_j)$, $j \in \{1, \ldots, N\}$, is not known at a central computational unit, but its knowledge is distributed among several nodes. An example of this setup can be the computation of a classifier for spam filtering [29], where the datum/label pairs are collected by the personal computers of $n$ users, and the $n$ computers may interact for computing the classifiers. For our numerical experiments we considered a problem in which features with label '+1' are sampled from the normal distribution with mean $10 \times \mathbf{1}_p$, while features with label '−1' are sampled from the normal distribution with mean $-10 \times \mathbf{1}_p$. After "sampling" the random constraints we distribute them among $n$ nodes. Then, we study the distributed computation of the solution to problem (18) on two network topologies: geometric random graph, and chain graph. The performance of ACC and VCC algorithms for $p = 4$ and $N = 20000$ total constraints is shown in Table III. The values shown in the table are the worst-case values over 20 runs of the algorithms. It can be seen that the number of iterations required for convergence of the ACC algorithm are linear in graph diameter, while they are equal to the graph diameter for the VCC algorithm. The number of constraints exchanged at each iteration among the nodes is small for ACC algorithm while this number is large for VCC algorithm.

*Example 5 (**Parallel linear classification**):* For the same set of data as in Example 4, we study the parallel computation of the optimal separating hyperplane. The parallel computation setup is modelled via a complete graph. The computation time of the ACC algorithm for parallel computation of the optimal separating hyperplane is shown in Fig. 6. The computation time is averaged over 20 runs of the algorithm. The computation time is shown, respectively, as a function of number of processors for $p = 4$ and $N = 200000$ total constraints, as a function of total number of constraints for $p = 4$ and $n = 50$ processors, and as a function of dimension $p$ for $N = 200000$ total constraints and $n = 50$ processors. In the first case, the minimum, average, and maximum number of active constraints are $2, 3.3$, and $5$, respectively, while the minimum, average, and maximum number of iterations are $4, 4.04$, and $5$, respectively. In the second case, the minimum, average, and maximum number of active constraints are $2, 3.09$, and $5$, respectively,

| | No. Nodes | Diameter | ACC | | VCC | |
|---|---|---|---|---|---|---|
| | | | Iter. | Constr. | Iter. | Constr. |
| Geometric random graph | 10 | 1 | 5 | | 1 | |
| | 50 | 2 | 11 | | 2 | |
| | 100 | 3 | 11 | 5 | 3 | 342 |
| | 500 | 5 | 24 | | 5 | |
| Chain graph | 10 | 10 | 37 | | 10 | |
| | 50 | 50 | 177 | | 50 | |
| | 100 | 100 | 319 | 5 | 100 | 365 |
| | 500 | 500 | 1498 | | 500 | |

TABLE III

DISTRIBUTED LINEAR CLASSIFICATION: MAXIMUM NUMBER OF ITERATIONS, MAXIMUM NUMBER OF EXCHANGED CONSTRAINTS, AND DIAMETER FOR DIFFERENT GRAPH TOPOLOGIES, FOR ACC AND VCC ALGORITHMS.

while the minimum, average, and maximum number of iterations are $4, 4.03$, and $6$, respectively. In the third case, the minimum, average, and maximum number of iterations are $4, 4.04$, and $5$, respectively, and the statistics of the constraints are shown in Fig. 6. It can be seen that the parallel computation of the optimal solution via ACC algorithm remarkably improves the computation time over the centralized computation. For large number of constraints, this improvement is significant even for a small number of processors.

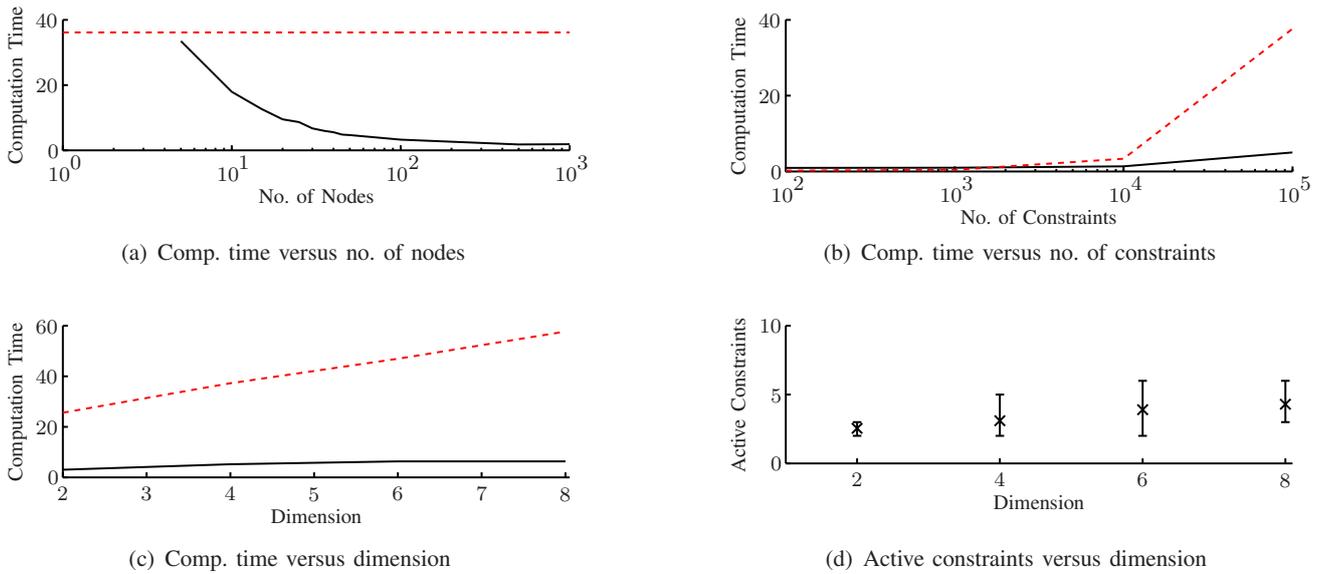

(a) Comp. time versus no. of nodes

(b) Comp. time versus no. of constraints

(c) Comp. time versus dimension

(d) Active constraints versus dimension

Fig. 6. Performance of ACC algorithm for parallel computation of the solution of linear classification problem. The solid black and dashed red lines represent parallel and centralized average computation time, respectively. The cross represents the average number of active constraints and the error bars represent the minimum and maximum active constraints.

## C. Parallel random model predictive control

Consider the LTI system

$$x_{t+1} = F(\xi)x_t + G(\xi)u_t + G_\gamma(\xi)\gamma_t, \tag{20}$$

where $t \in \mathbb{Z}_{\geq 0}$ is discrete time variable, $x_t \in \mathbb{R}^p$ is the system state, $u_t \in \mathbb{R}^q$ is the control input, $\gamma_t \in \Gamma \subset \mathbb{R}^{q_\gamma}$ is an unmeasured disturbance vector, $\xi \in \Xi \subseteq \mathbb{R}^w$ is vector of uncertain parameters, and $F(\xi) \in \mathbb{R}^{p \times p}, G(\xi) \in \mathbb{R}^{p \times q}, G_\gamma(\xi) \in \mathbb{R}^{p \times q_\gamma}$ are uncertain matrices. The design problem is to determine a control law that regulates the system state to some desired set, subject to some constraints on states and controls. In random model predictive control [30], one picks a control law of the form $u_t = K_f x_t + v_t$, where $K_f \in \mathbb{R}^{q \times p}$ is the static linear terminal controller gain and $v_t \in \mathbb{R}^q$ is the design variable. The design variable $v_t$ is picked to provide robustness with high probability. To determine the design variable that achieves such robustness, at each time $t$ and for a given finite horizon length $M$, $N$ realizations of the uncertain parameter $\xi$ and disturbance vectors $(\gamma_t, \ldots, \gamma_{t+M-1})$ are sampled and an optimization problem is solved. Let us denote these realizations by $(\xi^{(k)}, \gamma_t^{(k)}, \ldots, \gamma_{t+M-1}^{(k)}), k \in \{1, \ldots, N\}$, and define $g_t^{(k)} = [\gamma_t^{(k)\top} \ \ldots \ \gamma_{t+M-1}^{(k)\top}]^\top$, for each $k \in \{1, \ldots, N\}$. The design variable $v_t$ is determined by the solution of the following optimization problem:

$$
\begin{aligned}
\min_{\mathcal{V}_t} \quad & \max_{k \in \{1, \ldots, N\}} J(x_t, \mathcal{V}_t, \xi^{(k)}, g_t^{(k)}) \qquad \text{subject to :} \\
& f_X(x_{t+j|t}^{(k)}) \leq 0, \\
& f_U(K_f x_{t+j|t}^{(k)} + v_{t+j-1}) \leq 0, \\
& f_{X_f}(x_{t+M|t}^{(k)}) \leq 0, \\
& \text{for each } j \in \{1, \ldots, M\}, \text{ and for each } k \in \{1, \ldots, N\},
\end{aligned}
\tag{21}
$$

where $J : \mathbb{R}^p \times \mathbb{R}^{qM} \times \mathbb{R}^w \times \mathbb{R}^{q_\gamma M} \to \mathbb{R}$ is a cost function that is convex in $x_t$ and $\mathcal{V}_t$, $f_X : \mathbb{R}^p \to \mathbb{R}$, $f_U : \mathbb{R}^q \to \mathbb{R}$, and $f_{X_f} : \mathbb{R}^p \to \mathbb{R}$ are convex functions that capture constraints on the state at each time, the control at each time, and the final state, respectively, and

$$
\begin{aligned}
x_{t+j|t}^{(k)} &= (F_{\mathrm{cl}}(\xi^{(k)}))^{j-1} x_t + \Psi_j^{(k)} \mathcal{V}_t + \Upsilon_j^{(k)} g_t^{(k)} \\
F_{\mathrm{cl}}(\xi^{(k)}) &= F(\xi^{(k)}) + G(\xi^{(k)}) K_f \\
\Psi_j^{(k)} &= \left[ (F_{\mathrm{cl}}(\xi^{(k)}))^{j-1} G(\xi^{(k)}) \ \ldots \ F_{\mathrm{cl}}(\xi^{(k)}) G(\xi^{(k)}) \quad G(\xi^{(k)}) \quad 0 \ldots 0 \right] \in \mathbb{R}^{p \times qM} \\
\Upsilon_j^{(k)} &= \left[ (F_{\mathrm{cl}}(\xi^{(k)}))^{j-1} G_\gamma(\xi^{(k)}) \ \ldots \ F_{\mathrm{cl}}(\xi^{(k)}) G_\gamma(\xi^{(k)}) \quad G_\gamma(\xi^{(k)}) \quad 0 \ldots 0 \right] \in \mathbb{R}^{p \times q_\gamma M} \\
\mathcal{V}_t &= \left[ v_t^\top, \ldots, v_{t+M-1}^\top \right]^\top.
\end{aligned}
$$

Problem (21) is a random convex program of dimension $d = qM + 1$. Moreover, assuming that the problem admits a unique optimal solution with probability one and for $N > qM + 1$, for any realization of the parameter and the disturbance vector, the constraints on the state and the control are satisfied with expected probability at least $(N - qM)/(N + 1)$ [30]. Problem (21) is directly amenable to distributed solution via ACC algorithm. In the next section we consider the case in which the random constraints of the RCP are purposely distributed among $n$ processors that have to solve the problem in parallel fashion.

*1) Numerical results on parallel random MPC:* In order to achieve robustness with high probability, a large number of realizations of the parameter and disturbances are needed in the random convex program (21). This results in a large number of constraints and makes real-time centralized computation of the solution to the optimization problem (21) intractable. Therefore, we resort to the parallel computation of the solution the optimization problem (21) via ACC algorithm. We now apply the ACC algorithm to an example taken from [30], and show its effectiveness.

*Example 6 (Parallel random MPC):* Consider the LTI system (20) with

$$
F(\xi) = \begin{bmatrix} 1 + \xi_1 & \frac{1}{1+\xi_1} \\ 0.1\sin(\xi_4) & 1 + \xi_2 \end{bmatrix}, \quad G(\xi) = \begin{bmatrix} 0.3\arctan(\xi_5) \\ \frac{1}{1+\xi_3} \end{bmatrix}, \quad G_\gamma = \begin{bmatrix} 1 & 0 \\ 0 & 1 \end{bmatrix},
$$

where each of the random parameters $\xi_1, \xi_2, \xi_3$ is uniformly distributed in the interval $[-0.1, 0.1]$, while $\xi_4, \xi_5$ are distributed according to Gaussian distributions with zero mean and unit variance. Let the horizon be $M = 10$ and the uncertainty $\gamma$ be uniformly distributed over set $\Gamma = \{\gamma \in \mathbb{R}^2 : \|\gamma\|_\infty\} \leq 0.05$. Assume that $f_X(x) = \|x\|_\infty - 10$, $f_U(u) = |u| - 5$, and $f_{X_f}(z) = \|z\|_\infty - 1$. Given the terminal controller gain $K_f = [-0.72 \ -1.70]$ and the cost function $J(x_t, \mathcal{V}_t, \xi^{(k)}, g_t^{(k)}) = \max_{j \in \{1, \ldots, M\}} \max\{0, \|x_{t+j-1|t}^{(k)}\|_\infty - 1\} + \|\mathcal{V}_t\|_2^2$. For this set of data, the computation time of the ACC algorithm averaged over 20 runs of the algorithm for parallel computation of the solution to optimization problem (21) is shown in Fig. 7. The computation time is shown, respectively, as a function of number of processors for 1000 realizations of the random parameters, and as a function of number of realizations of the random parameters for 50 processors. In the first case, the minimum, average, and maximum number of active constraints are $2, 2.55,$ and $6$, respectively, while the minimum, average, and maximum number of iterations are $3, 3.73,$ and $5$, respectively. In the second case, the minimum, average, and maximum number of active constraints are $2, 2.18,$ and $4$, respectively, while the minimum, average, and maximum number of iterations are $3, 4.03,$ and $5$, respectively.

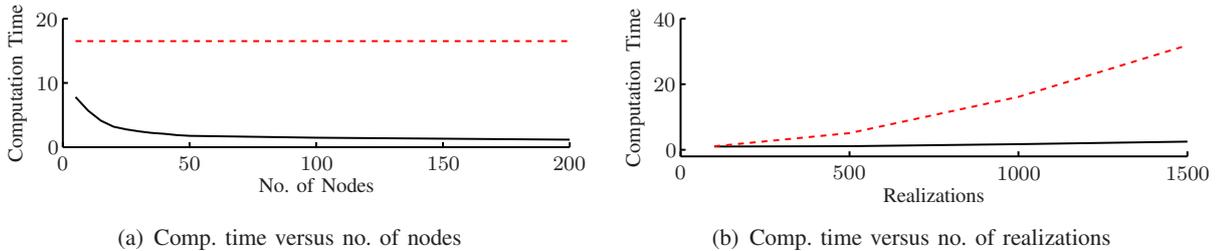

(a) Comp. time versus no. of nodes

(b) Comp. time versus no. of realizations

Fig. 7. Performance of the ACC algorithm for parallel random model predictive control. The solid black and dashed red lines represent parallel and centralized average computation time, respectively.

## D. Example of distributed outliers rejection

We conclude the numerical part of this paper with a brief example of distributed constraints removal, applied to the distributed estimation problem presented in Section VII-A1. We consider $n = 50$ sensors measuring a random variable $\theta$, using the same measurement model of Example 1 (homogeneous sensor network). The overall number of measurements (acquired by all nodes) is $N = 3000$. The original scenario solution that satisfies all $N = 3000$ constraints can assure a violation probability smaller than $\epsilon = 10^{-2}$ with confidence level greater than $1 - \beta = 1 - 2 \times 10^{-8}$. According to RCP theory we can remove $r = 165$ constraints, still guaranteeing that the violation probability is smaller that $10^{-1}$ with confidence level $1 - \beta$ close to $1$. Therefore the nodes apply Algorithm 4 (the ACC algorithm is used within the removal strategy), computing a scenario solution which satisfies all but $r = 165$ constraints. Thus, with a little compromise over the bound on the violation probability, the constraints removal allows reducing the size of the ellipsoid, hence improving the informativeness of the confidence ellipsoid. In Fig. 8, we report the confidence ellipsoids computed at one node using Algorithm 4, after rejecting number of outliers $\eta = \{0, 20, 40, \ldots, 140, 160\}$, together with the final ellipsoid satisfying all but $r = 165$ constraints.

## VIII. CONCLUSION

In this paper, we studied distributed computation of the solution to random convex program (RCP) instances. We considered the case in which each node of a network of processors has local knowledge of only a subset of constraints of the RCP, and the nodes cooperate in order to reach the solution of the *global* problem (i.e., the problem including all constraints). We proposed two distributed algorithms, namely, the active constraints consensus (ACC) algorithm and vertex constraints consensus (VCC) algorithm. The ACC algorithm computes the solution in finite time and requires the nodes to exchange a small number of constraints at each iteration. Moreover, a parallel implementation of the ACC algorithm remarkably improves the computational effort compared to a centralized solution of the RCP. The VCC algorithm converges to the solution in a number of iterations equal to the graph diameter. We also developed a variant of VCC algorithm, namely, quantized vertex constraints consensus (qVCC), that restricts the number of constraints to be exchanged at each iteration. We further proposed a distributed constraints removal strategy for outlier rejection within the framework of RCP with violated constraints. Finally, we presented several applications of the proposed distributed algorithms, including estimation, classification, and random model predictive control.

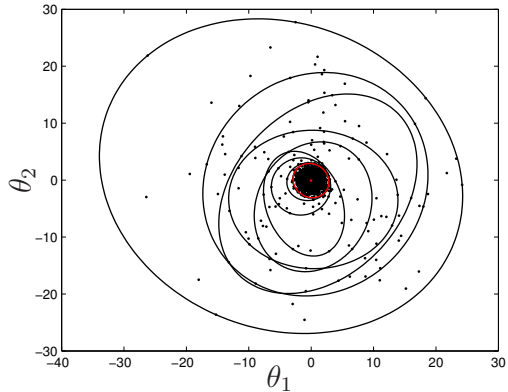

Fig. 8. Measurements taken by all the sensors in the network (black dots) and confidence ellipsoids at one node after rejecting number of outliers $\eta = \{0, 20, 40, \ldots, 140, 160\}$ in Algorithm 4. The red ellipsoid is the one produced after discarding $r = 165$ measurements according to the *distributed constraints removal procedure*.

## APPENDIX

### A.1: Proof of Proposition 2

We start by establishing the first statement. Let $c$ be a support constraint for a feasible problem in the form (1). Call $\hat{x}^* = x^*(C)$ and $\check{x}^* = x^*(C\backslash\{c\})$. From the definition of support constraints, it follows that $a^\top \check{x}^* < a^\top \hat{x}^*$. Assume by contradiction that $c$ is not active at $\hat{x}^*$, i.e. that $f_c(\hat{x}^*) < 0$. Consider a point $\bar{x}$ on the segment connecting $\hat{x}^*$ and $\check{x}^*$: $\bar{x}(\lambda) = \lambda\check{x}^* + (1-\lambda)\hat{x}^*$, $\lambda \in [0, 1]$. It follows immediately that $a^\top \bar{x}(\lambda) < a^\top \hat{x}^*$, for every $\lambda \in [0, 1)$. Moreover, by convexity, $\bar{x}(\lambda)$ satisfies all constraints, except possibly constraint $c$. However, since $\hat{x}^*$ is in the interior of the convex set defined by $f_c \leq 0$, there must exist values of $\lambda$ sufficiently small such that $\bar{x}(\lambda)$ satisfies also $f_c(\bar{x}(\lambda)) \leq 0$. But then $\bar{x}(\lambda)$ would satisfy all constraints and yield an objective value that improves upon that of $\hat{x}^*$. This contradicts optimality of $\hat{x}^*$ and hence proves that $c$ must be active at $\hat{x}^*$.

We now establish the second statement. We first demonstrate that each essential set $E_i \doteq \text{Es}_i(C)$ needs be irreducible, i.e., $E_i = \text{Sc}(E_i)$. By definition, each $E_i$ is a minimum cardinality set satisfying $J^*(E_i) = J^*(C)$. Now assume by contradiction that there exists a constraint $c \in E_i$, such that $J^*(E_i) = J^*(E_i\backslash\{c\})$. This implies that there exists a set $E_i\backslash\{c\}$, which is also invariant for $C$, i.e., $J^*(E_i\backslash\{c\}) = J^*(E_i) = J^*(C)$, and has smaller cardinality than $E_i$, leading to contradiction. Now we can prove the statement: if each constraint in $\text{Es}_i(C)$ is a support constraint for problem $P[\text{Es}_i(C)]$, it needs to be active for the problem $P[\text{Es}_i(C)]$, see claim (i). Consequently, if $x_i^*$ is the optimal solution for $P[\text{Es}_i(C)]$, then $f_j(x_i^*) = 0$, $\forall j \in \text{Es}_i(C)$. From the unique minimum condition and locality, it follows that

$$J^*(\text{Es}_i(C)) = J^*(C) \implies x^*(\text{Es}_i(C)) = x^*(C),$$

for each $i \in \{1, \ldots, n_e\}$. Therefore, $f_j(x^*(C)) = 0$, for each $j \in \text{Es}_i(C)$, $i \in \{1, \ldots, n_e\}$, and $\text{Ac}(C) \supseteq \cup_{i=1}^{n_e} \text{Es}_i(C)$. $\qquad\square$

### A.2: Proof of Proposition 3

We start by establishing the first statement. According to the update rule of the ACC algorithm, the sequence of local optimal objective $J_i^*(t)$ satisfies

$$
\begin{aligned}
J_i^*(t+1) \doteq J^*\big(L_i(t+1)\big) &= J^*\big(A_i(t) \cup (\cup_{j \in \mathcal{N}_{\text{in}}(i)} A_j(t)) \cup C_i\big) \\
\text{[by monotonicity]} &\geq J^*\big(A_i(t)\big) \\
\text{[by Corollary 1]} &= J^*\big(L_i(t)\big) = J_i^*(t),
\end{aligned}
$$

then $J_i^*(t)$ is non-decreasing in $t$.

The proof of the second statement is more involved and works as follows. We first observe that, for any directed edge $(i, j)$ it holds

$$
\begin{aligned}
J_j^*(t+1) \doteq J^*\big(L_j(t+1)\big) &= J^*\big(A_j(t) \cup (\cup_{k \in \mathcal{N}_{\text{in}}(j)} A_k(t)) \cup C_j\big) \\
\text{[by monotonicity and } i \in \mathcal{N}_{\text{in}}(j)] &\geq J^*\big(A_i(t)\big) \\
\text{[by Corollary 1]} &= J^*\big(L_i(t)\big) = J_i^*(t),
\end{aligned}
$$

which can be easily generalized to a generic pair of nodes $i$, $j$ connected by a directed path of length $l_{ij}$ (such path always exists for the hypothesis of strong connectivity):

$$
J_j^*(t+l_{ij}) \geq J_i^*(t). \tag{22}
$$

Moreover, we demonstrate that for any directed edge $(i, j)$ it holds

$$
J_j^*(t+1) = J_i^*(t) \quad \Longleftrightarrow \quad x_j^*(t+1) = x_i^*(t). \tag{23}
$$

The reverse implication in (23) is straightforward, since the objective function is the same for both nodes. The direct implication is again trivial in the infeasible case, while for $J_j^*(t+1) = J_i^*(t) < \infty$ it can be proven as follows. For the uniqueness condition, adding a constraint $c$ that is not satisfied at (or *violates*) $x_j^*(t+1)$ leads to an increase in the objective value, i.e., $J^*(L_j(t+1) \cup \{c\}) > J^*(L_j(t+1))$. Now, since $L_j(t+1) \supseteq A_i(t)$, and $J^*(L_j(t+1)) = J_j^*(t+1) = J_i^*(t) = J^*(A_i(t))$, by locality, if $J^*(L_j(t+1) \cup \{c\}) > J^*(L_j(t+1))$, then $J^*(A_i(t) \cup \{c\}) > J^*(A_i(t))$, which implies that also $x_i^*(t)$ is violated by $c$. Therefore, we concluded that every constraint that violates $x_j^*(t+1)$ also violates $x_i^*(t)$ and this may happen only if $x_j^*(t+1) = x_i^*(t)$. Again the correspondence between objective values and optimal solutions can be easily generalized to a generic pair of nodes $i$, $j$ connected by a directed path of length $l_{ij}$:

$$
J_j^*(t+l_{ij}) = J_i^*(t) \quad \Longleftrightarrow \quad x_j^*(t+l_{ij}) = x_i^*(t). \tag{24}
$$

We now claim that the objective at one node cannot remain the same for $2\texttt{diam}(\mathcal{G})+1$ iterations, unless the algorithm has converged. In the infeasible case the proof is trivial: according to the update rule of the ACC if node $i$ has detected an infeasible local problem, i.e., $J_i^*(t) = \infty$, it directly stops the execution of the algorithm since it is already sure of detaining the global solution. Let us instead consider the feasible case. We assume by contradiction that $J_i^*(t) = J_i^*(t+2\texttt{diam}(\mathcal{G})) < \infty$ and there exists a node $j$ with at least a constraint that is not satisfied by $x_i^*(t) = x_i^*(t+2\texttt{diam}(\mathcal{G}))$. Let us consider a directed path of length $l_{ij}$ from $i$ to $j$: we already observed in (22) that $J_j^*(t+l_{ij}) \geq J_i^*(t)$. However, since there are constraints at node $j$ that violates $x_i^*(t)$, equality cannot hold, see (24), and $J_j^*(t+l_{ij}) > J_i^*(t)$. By definition, the length $l_{ij}$ of the path from $i$ to $j$ is bounded by graph diameter and the local objective is non-decreasing, therefore $J_j^*(t+\texttt{diam}(\mathcal{G})) > J_i^*(t)$. Now consider the path from $j$ to $i$ of length $l_{ji}$: according to (22) it must hold $J_j^*(t+\texttt{diam}(\mathcal{G})) \leq J_i^*(t+\texttt{diam}(\mathcal{G})+l_{ji}) \leq J_i^*(t+2\texttt{diam}(\mathcal{G}))$. Using the two inequalities found so far we obtain $J_i^*(t) < J_j^*(t+\texttt{diam}(\mathcal{G})) \leq J_i^*(t+2\texttt{diam}(\mathcal{G}))$, which contradicts the assumption that the objective at node $i$ remains constant for $2\texttt{diam}(\mathcal{G})+1$ iterations. Therefore, before convergence the local objective $J_i^*(t)$ has to be strictly increasing every $2\texttt{diam}(\mathcal{G})+1$ iterations. Moreover, the sequence $J_i^*(t)$ is upper bounded, since, by monotonicity, for any $L \subseteq C$, $J^*(L) \leq J^*(C)$, and $J_i^*(t)$ can assume a finite number of values, i.e., $J^* \in \mathcal{J} \doteq \{J^*(L) : L \subseteq C\}$; therefore the sequence converges to a constant value, say $J_i^*(T)$, in finite time. We now demonstrate that after convergence all nodes need to have the same local objective, i.e., $J_i^*(T) = \hat{J}$, for each $i \in \{1, \ldots, n\}$. For simplicity of notation, we drop the time index in the following discussion. Assume by contradiction that two nodes, say $i$ and $j$, have different objective values, $J_i^* > J_j^*$. From the assumption of strongly connectivity of the graph $\mathcal{G}$, there exist a directed path between $i$ and $j$. Using relation (22) we obtain $J_i^* \leq J_j^*$, leading to contradiction. Therefore, for any pair of nodes $i$ and $j$, it must hold that $J_i^* = J_j^* = \hat{J}$, implying $J_i^* = \hat{J}$, for each $i \in \{1, \ldots, n\}$. With a similar reasoning, and using (24), we can also conclude that $J_i^* = \hat{J}$, for each $i \in \{1, \dots, n\}$, implies $x_i^* = \hat{x}$, for each $i \in \{1, \dots, n\}$. Now it remains to show that the local objectives $\hat{J}$ and the local solutions $\hat{x}$ actually coincide with $J^*(C)$ and $x^*(C)$. In the infeasible case this is again trivial: if the local objectives coincide with $\hat{J} = \infty$, by monotonicity the global problem cannot be either than infeasible, then $J^*(C) = \hat{J} = \infty$ and $x^*(C) = \hat{x} = \texttt{NaN}$. The feasible case can be proven as follows. If all nodes have the same local solution $\hat{x}$, it means that (i) $\hat{x}$ satisfies the local constraint set $C_i$, $i \in \{1, \dots, n\}$, which implies that $\hat{x}$ is feasible for the global problem. Moreover, by monotonicity, $\hat{J} \le J^*(C)$ (since $\hat{J}$ is the optimal value of a subproblem having constraint set $L \subseteq C$). Assume by contradiction that $\hat{J} < J^*(C)$, which implies that (ii) $\hat{J} = a^\top \hat{x} < a^\top x^*(C) = J^*(C)$; therefore $\hat{x}$ attains a smaller objective than $x^*(C)$, see (ii), and satisfies all constraints in $C$, see (i), contradicting the optimality of $x^*(C)$. Therefore it must hold $\hat{J} = J^*(C)$. With the same reasoning we used for proving (23), we also conclude that $\hat{x} = x^*(C)$.

To prove the third statement, we show that the set $A_i$ contains all the constraints that are globally active for $P[C]$. If $J_i^* = J^*(C) = \infty$ the implication is trivial, since $A_i = \texttt{Ac}(C) = \emptyset$. In the feasible case the proof proceeds as follows. According to the second statement, we have $x_i^* = x^*(A_i) = x^*(C)$, $i \in \{1, \dots, n\}$. By contradiction, let us suppose that there exists a globally active constraint $c$ that is contained in the local constraint set $C_i$ of a node $i$, but is not in the candidate set $A_j$ of node $j$. Let us consider a directed path from $i$ to $j$ and relabel the nodes in this path from $1$ to $l$. Starting from node $1$ we observe that, since $x_1^* = x^*(C)$ and $c$ is active for $P[C]$, then $c \in A_1$. At each iteration of the active constraint consensus, node $2$ in the path computes $A_2 = \texttt{Ac}(A_2 \cup (\cup_{j \in \mathcal{N}_{\mathrm{in}}(2,t)} A_j) \cup C_2)$. Therefore, since $c \in A_1$ and $x_1^* = x_2^*$, it holds $c \in A_2$. Iterating this reasoning along the path from $i$ to $j$ we conclude that $c \in A_j$ leading to contradiction.

To prove the fourth statement, we observe that, if the local problem at node $i$ is infeasible, then the node only has to transmit its local objective, $J_i(t)^* = \infty$, since the candidate set $A_i(t)$ is empty. If the local problem $P[L_i]$ is feasible, then the unique minimum condition assures that the minimum is attained at a single point, say $x^*(L_i)$. If constraints are in general position, then no more than $d$ constraints may be tight at $x^*(L_i)$, hence at most $d$ constraints are active. Therefore, in the feasible case, the number of constraints to be transmitted is upper-bounded by $d$. □

### A.3: Proof of Proposition 4

We start by recalling a basic property which is a direct consequence of the definition of the feasible set: for any set of constraints $C_1$ and $C_2$, it holds:

$$\texttt{Sat}(C_1) \cap \texttt{Sat}(C_2) = \texttt{Sat}(C_1 \cup C_2). \tag{25}$$

To prove the first statement, we consider a generic node $i$. At time $t$ node $i$ receives the candidate sets from the incoming neighbors and computes $V_i(t+1) = \texttt{vert}\big(L_i(t+1)\big) = \texttt{vert}\big(V_i(t) \cup \big(\cup_{j \in \mathcal{N}_{\mathrm{in}}(i)} V_j(t)\big)\big)$. It follows that

$$
\begin{aligned}
\texttt{Sat}(V_i(t+1)) &= \texttt{Sat}\big(\texttt{vert}\big(V_i(t) \cup \big(\cup_{j \in \mathcal{N}_{in}(i)} V_j(t)\big)\big)\big) \\
\text{[by Fact 1]} &= \texttt{Sat}\big(V_i(t) \cup \big(\cup_{j \in \mathcal{N}_{in}(i,t)} V_j(t)\big)\big) \\
\text{[by equation (25)]} &= \texttt{Sat}(V_i(t)) \cap \big(\cap_{j \in \mathcal{N}_{in}(i)} \texttt{Sat}(V_j(t))\big) \subseteq \texttt{Sat}(V_i(t)).
\end{aligned}
\tag{26}
$$

If $\texttt{Sat}(V_i(t)) = \emptyset$ (infeasible local problem) also $\texttt{Sat}(V_i(t+1)) = \emptyset$, according to (26), then $J_i^*(t+1) = J_i^*(t) = \infty$, and the objective is non-decreasing. If $\texttt{Sat}(V_i(t)) \ne \emptyset$ (feasible local problem) we can prove the statement as follows. Assume by contradiction that there exists $\bar{x} \in \texttt{Sat}(V_i(t+1))$ such that $a^\top \bar{x} \doteq J^*(V_i(t+1)) < J^*(V_i(t))$. Equation (26) assures that $\texttt{Sat}(V_i(t+1)) \subseteq \texttt{Sat}(V_i(t))$, therefore, $\bar{x} \in \texttt{Sat}(V_i(t))$ and exists a point in the feasible set of problem $P[V_i(t)]$, whose value is smaller than $J^*(V_i(t))$. This contradicts the optimality of $J^*(V_i(t))$. Therefore, it must hold $J^*(V_i(t+1)) \ge J^*(V_i(t))$.

To prove the second statement, we show that after $T \doteq \texttt{diam}(\mathcal{G})$ iterations a generic node $i$ satisfies $\texttt{Sat}(V_i(T)) = \texttt{Sat}(C)$. Consider a generic node $j$ and a directed path from a node $j$ to node $i$ (this path does exist for the hypothesis of strong connectivity). We relabel the nodes on this path from 1 to $l$, such that the last node is $i$. Node 1 initializes $V_1(0) = \mathtt{vert}(C_1)$, then $\mathtt{Sat}(V_1(0)) = \mathtt{Sat}(C_1)$. At the first iteration, node 2 computes $V_2(1) = \mathtt{vert}\big(V_2(0) \cup (\cup_{j \in \mathcal{N}_{\mathrm{in}}(2)} V_j(0))\big)$. Since node 1 is in $\mathcal{N}_{\mathrm{in}}(2)$, it follows from equation (26) that $\mathtt{Sat}(V_2(1)) \subseteq \mathtt{Sat}(V_1(0))$. Repeating the same reasoning along the path, and for the original labeling, we can easily prove that $\mathtt{Sat}(V_i(l_{ij})) \subseteq \mathtt{Sat}(V_j(0)) = \mathtt{Sat}(C_j)$, where $l_{ij}$ is the distance between $i$ and $j$. Therefore, after a number of iterations equal to the distance between $j$ and $i$, every feasible solution at node $i$ satisfies the constraints of node $j$. Since the maximum distance between $i$ and any other node is the diameter of the graph, in $T \doteq \mathtt{diam}(\mathcal{G})$ iterations, node $i$ satisfies $\mathtt{Sat}(V_i(T)) \subseteq \mathtt{Sat}(C_j)$ for all $j$. Since this last property holds for all $j$, it also holds $\mathtt{Sat}(V_i(T)) \subseteq \cap_{j \in \{1,\dots,n\}} \mathtt{Sat}(C_j) = \mathtt{Sat}(C)$. However, $V_i(T)$ is a subset of $C$, and it follows that $\mathtt{Sat}(V_i(T)) \supseteq \mathtt{Sat}(C)$. Thus, $\mathtt{Sat}(V_i(T)) = \mathtt{Sat}(C)$. Since the local problem $P[V_i(T)]$ and the global problem $P[C]$ have the same objective direction and the same feasible set they attain the same (unique) solution, i.e., $x^*(V_i(T)) = x^*(C)$.

We now establish the third statement. We note that $V_i(T) = \mathtt{vert}(C)$ is a direct consequence of the update rule of the VCC algorithm. To prove the latter part of the statement, we assume by contradiction that $c \in C$ is a support constraint for $P[C]$ but $c \notin \mathtt{vert}(C)$. The relation $c \notin \mathtt{vert}(C)$ implies that $\mathtt{vert}(C) \subseteq C \backslash \{c\}$. It follows from monotonicity that (i) $J^*(\mathtt{vert}(C)) \leq J^*(C \backslash \{c\})$. According to Fact 1 it also holds (ii) $J^*(\mathtt{vert}(C)) = J^*(C)$. Combining (i) and (ii), we obtain $J^*(C) \leq J^*(C \backslash \{c\})$. By monotonicity, it cannot be $J^*(C) < J^*(C \backslash \{c\})$, then $J^*(C) = J^*(C \backslash \{c\})$, but this contradicts the assumption that $c$ is a support constraint. $\qquad \square$

### A.4: Proof of Proposition 5

The proof of the first and the third statement follows similar to the proof of the first and third statement in Proposition 3. We now establish the second statement. Similar to the VCC algorithm, we show that after $T \leq \sum_{k=0}^{\mathtt{diam}(\mathcal{G})-1} \lceil \frac{N_{\max}(d_{\max}+1)^k}{m} \rceil$ iterations a generic node $i$ satisfies $\mathtt{Sat}(V_i(T)) = \mathtt{Sat}(C)$. Consider a generic pair of nodes $i$, $j$ and a directed path of length $l_{ji}$ from $j$ to $i$ (this path does exist for the hypothesis of strong connectivity). Relabel the nodes on this path from 1 to $l$, such that the last node is $i$. We observe that, after the initialization, the local candidate set $V_1(0) = T_1(0) = \mathtt{vert}(C_1)$ has cardinality $|T_1(0)| \leq N_{\max}$. Since the transmission set is managed using a FIFO policy, after at most $\lceil \frac{N_{\max}}{m} \rceil$ communication rounds the node has transmitted all the constraints in $V_1(0)$ to node 2. Therefore, $\mathtt{Sat}(V_2(\lceil \frac{N_{\max}}{m} \rceil)) \subseteq \mathtt{Sat}(V_1(0)) = \mathtt{Sat}(C_1)$. Moreover, $\left| V_2(\lceil \frac{N_{\max}}{m} \rceil) \right| \leq \sum_{j \in \mathcal{N}_{\mathrm{in}}(2) \cup \{2\}} N_j \leq N_{\max}(d_{\max}+1)$ (worst case in which the incoming neighbors have to transmit all their local constraints and all constraints are vertices of the convex hull). After at most $\lceil \frac{N_{\max}(d_{\max}+1)}{m} \rceil$ further iterations, node 2 has transmitted all constraints in $V_2(\lceil \frac{N_{\max}}{m} \rceil)$ to node 3. Therefore, $\mathtt{Sat}(V_3(\lceil \frac{N_{\max}}{m} \rceil + \lceil \frac{N_{\max}(d_{\max}+1)}{m} \rceil)) \subseteq \mathtt{Sat}(V_2(\lceil \frac{N_{\max}}{m} \rceil)) \subseteq \mathtt{Sat}(C_1)$. Also, $\left| V_3(\lceil \frac{N_{\max}}{m} \rceil + \lceil \frac{N_{\max}(d_{\max}+1)}{m} \rceil) \right| \leq \sum_{j \in \mathcal{N}_{\mathrm{in}}(3) \cup \{3\}} \left| V_j(\lceil \frac{N_{\max}}{m} \rceil) \right| \leq N_{\max}(d_{\max}+1)^2$. Repeating the same reasoning along the directed path, for the original labeling, we obtain $\mathtt{Sat}\left(V_i\big(\sum_{k=0}^{l_{ji}-1} \lceil \frac{N_{\max}(d_{\max}+1)^k}{m} \rceil\big)\right) \subseteq \mathtt{Sat}(C_j)$. Therefore, every feasible solution at node $i$ satisfies the constraints of node $j$ at distance $l_{ji}$ in a number of iterations no larger than $\sum_{k=0}^{l_{ji}-1} \lceil \frac{N_{\max}(d_{\max}+1)^k}{m} \rceil$. Since the maximum distance between $i$ and any other node is the diameter of the graph, it follows that in $T \leq \sum_{k=0}^{\mathtt{diam}(\mathcal{G})-1} \lceil \frac{N_{\max}(d_{\max}+1)^k}{m} \rceil$ iterations node $i$ satisfies $\mathtt{Sat}(V_i(T)) \subseteq \mathtt{Sat}(C_j)$ for all $j$. Since this property holds for all $j$, it also holds $\mathtt{Sat}(V_i(T)) \subseteq \cap_{j \in \{1,\dots,n\}} \mathtt{Sat}(C_j) = \mathtt{Sat}(C)$. Since $V_i(T)$ is a subset of $C$, $\mathtt{Sat}(V_i(T)) \supseteq \mathtt{Sat}(C)$. Therefore, $\mathtt{Sat}(V_i(T)) = \mathtt{Sat}(C)$. Finally, $T$ can be rewritten as

$$\sum_{k=0}^{\mathtt{diam}(\mathcal{G})-1} \left\lceil \frac{N_{\max}(d_{\max}+1)^k}{m} \right\rceil \leq \left\lceil \frac{N_{\max}}{m} \right\rceil \sum_{i=0}^{\mathtt{diam}(\mathcal{G})-1} (d_{\max}+1)^k =$$

$$\left\lceil \frac{N_{\max}}{m} \right\rceil \frac{1 - (d_{\max}+1)^{\mathtt{diam}(\mathcal{G})}}{1 - (d_{\max}+1)} = \left\lceil \frac{N_{\max}}{m} \right\rceil \frac{(d_{\max}+1)^{\mathtt{diam}(\mathcal{G})}-1}{d_{\max}},$$

which coincides with the bound in the second statement. Since the local problem $P[V_i(T)]$ and the global problem $P[C]$ have the same objective direction and the same feasible set they attain the same (unique) solution, i.e., $x^*(V_i(T)) = x^*(C)$. $\qquad\square$